\documentclass[11pt, a4paper, amsart]{article}

\usepackage{amsmath, amssymb, theorem, latexsym, bm}

\usepackage{mathptmx}  

\usepackage[english, francais]{babel}
\newtheorem{theorem}{Theorem}[section]
\newtheorem{lemma}[theorem]{Lemma}

\newtheorem{proposition}[theorem]{Proposition}

\newtheorem{corollary}[theorem]{Corollary}

\newtheorem{remark}[theorem]{Remark}

\allowdisplaybreaks
\renewcommand{\thefootnote}{\fnsymbol{footnote}}

\numberwithin{equation}{section}

\title {\textsf {On Positive Harmonic Functions in Cones and Cylinders
\\ }
}
\author { Alano Ancona{$^{}$}\\ {\small D\'{e}partement de Math\'{e}matiques, B\^atiment 425, Universit\'{e} Paris-Sud 11}
\\ {\small  Orsay 91\hspace{0.5pt}405 France}} 
\date{{\small \today}} 

\begin{document}

\selectlanguage{english}

\maketitle



{\renewcommand\thefootnote{}
\footnote{2010 Mathematics Subject Classification 31C35, 60J45, 60J60, 45C05
}}

\vspace{2truemm}
\selectlanguage{english}

\renewcommand{\abstractname}

\noindent{\bf Abstract.} {\footnotesize We first consider a question raised by Alexander Eremenko and show that if $\Omega $ is an arbitrary connected open   cone in ${ \mathbb  R}^d$, then any two positive harmonic functions in $\Omega $ that vanish on $\partial \Omega $ must be proportional -an already known fact when $\Omega $ has a Lipschitz basis or more generally a John basis.  It is also shown however that  when $d \geq 4$, there can be more than one Martin point  at infinity for the cone  though  non-tangential convergence to the canonical Martin point at infinity  always holds. In contrast,  when  $d \leq 3$, the Martin point at infinity is unique for every cone. These properties connected with the dimension are  related to well-known results of M.\ Cranston and T.\ R.\ McConnell about the  lifetime of conditioned Brownian motions in  planar domains and also to  subsequent results by R. Ba\~nuelos and B. \nolinebreak  Davis. We also investigate the nature of the Martin points arising at infinity as well as the effects on the Martin boundary resulting from the existence of John cuts in the basis of the cone or from other regularity assumptions. The main results together with their  proofs   extend  to  cylinders  ${ \mathcal C}_Y( \Sigma )= { \mathbb R} \times  \Sigma$ --where $ \Sigma $ is a relatively compact region of a manifold $M$--, equipped with a suitable  second order  elliptic operator.}

\footnote[0]{\small $^1$ {\em \!Email address}: {alano.ancona@math.u-psud.fr}}

\section{Introduction. Main results.}\label{sectintro}

We consider the cone 
 $ { \mathcal C}_o( \Sigma ) $  of ${ \mathbb R} ^d$, $d \geq 2$, generated by a region $ \Sigma $ of the unit sphere $S_{d-1}$, i.e., ${ \mathcal C}_o( \Sigma )= \{ r\omega \,;\, r>0,\; \omega   \in    \Sigma \; \}$, and study the positive  harmonic functions in ${ \mathcal C}_o( \Sigma )$ (where ${ \mathcal C}_o$ is for cone). 
 Recently Alexander Eremenko  asked   whether it is always true that any two such functions that moreover vanish on $\partial { \mathcal C}_o( \Sigma )$ must be  proportional. Our first main result, Theorem \ref{theorem1} below (see also Theorem \ref{theorem2bis}), answers this question by the positive. A generalization  to a large class of cylinders is described in section \ref{sectgene}.

To deal with   non necessarily Dirichlet-regular $ \Sigma $, we  say, following a usual convention, that a function  $w$   in $ \Sigma $ vanishes  on  the open subset $ T $ of $\partial  \Sigma $ (or, more precisely, that $w$ vanishes in the weak sense on $T$) if $w$ is bounded in a neighborhood of each $ \xi   \in    T   $ and if $ A:=$ $\{  \xi   \in    T  \,;\,  \displaystyle \limsup _{ \Sigma  \ni x \to   \xi  } \vert  w(x) \vert   >0\,  \}$ is polar  in $S_{d-1}$. By definition, $A \subset   S_{d-1}$ is polar in $S_{d-1}$ if  for each $  \xi     \in   A$ there is a chart of $S_{d-1}$, $\chi: V \to  W \subset   { \mathbb R} ^{d-1}$, $V\ni \xi  $,  such that $\chi (V \cap    A)$ is  polar in ${ \mathbb R} ^{d-1}$. Note that  $ A  $ is polar in $S_{d-1}$ if and only if $ \{ t \xi  \,;\, \xi   \in   A,\, t> 0 \} $ is polar in ${ \mathbb R} ^d$.  

For a function $w$ defined in 	a region $\Omega $  of ${ \mathbb R} ^d$,  the relation  $w=0$ on  $T \subset   \partial \Omega $ is defined similarly.
If $w:\Omega  \to  { \mathbb R} $  is harmonic with respect to a second order uniformly elliptic operator  in divergence form with bounded measurable  coefficients in  $\Omega $ and if $w=0$ in an open subset $W$ of $  \partial V$ then $ \displaystyle \lim_{x \to   \xi  }w(x)=0$  for every  Dirichlet-regular boundary point $  \xi     \in   W $.

\begin{theorem} \label{theorem1} The nonnegative harmonic functions  in ${ \mathcal C}_o( \Sigma )$	 which vanish (in the weak sense)   on the boundary of ${ \mathcal C}_o( \Sigma )$ are the functions $h$ in the form $h(r\omega )= c\, r^ {\alpha_{ \Sigma }} \varphi_0  (\omega )$, $\omega  \in    \Sigma $, $r>0$, where $c $ is a nonnegative constant,  $ \alpha _ \Sigma = {\frac {-(d-2)+\sqrt{(d-2)^2+4 \lambda _1( \Sigma )} } {2 } }$ and $ \varphi  _0$ is a  positive solution of $ \Delta  _{S_{d-1}}  \varphi  _0+ \lambda _1 ( \Sigma )\,\varphi  _0=0$ in $ \Sigma $.   
\end{theorem}

Here $ \Delta  _{S_{d-1}}$ is  the spherical  Laplacian --denoted also $ \Delta  _S$ in the rest of the paper-- and  $  \lambda _1( \Sigma )$ (later denoted $ \lambda _1$) is the first eigenvalue of the opposite of the Dirichlet Laplacian in $ \Sigma $.  As well-known $ \lambda _1( \Sigma )$ coincides with the Raleigh constant of $ \Sigma $, i.e.\  $ \lambda _1( \Sigma )=\inf  \{ \int  \vert \nabla u \vert  ^2\, d \sigma  _ {S_{d-1}} ;$ $u \in   C_c^1( \Sigma ),\,\, \int  \vert  u \vert  ^2\, d \sigma  _{S_{d-1}} \geq 1\,  \}$ --where $ \sigma  _ {S_{d-1}} $ is the standard Riemannian spherical measure in $S_{d-1}$-- and $ \lambda _1( \Sigma )$ is $>0$ if and only if  $S_{d-1}\setminus  \Sigma $ is not polar in $S_{d-1}$ (see e.g.\ \cite{GZ}).  In this case, $ {\frac {1 } {  \lambda _1( \Sigma ) } } $ is also the largest eigenvalue of the (nonnegative self-adjoint compact) Green's operator in $L^2( \Sigma ; \sigma  _ {S_{d-1}}) $, $ \varphi \mapsto G( \varphi  )=(-\Delta  _S)^{-1}( \varphi  )$. It is known (see e.g.\ \cite{ancbook}) that $ \lambda _1( \Sigma )$ is  also the greatest real $ \lambda $ for which there is a positive $ (\Delta  _S+ \lambda I)$-superharmonic function in $ \Sigma $ (distinct from the constant $+ \infty $). For $ \lambda = \lambda _1( \Sigma )$ such a function  is  unique --up to multiplication by a  constant-- and there is a unique  positive solution $ \varphi  _ 0\in   H_0^1( \Sigma )$ of $ \Delta  _S \varphi  _0+ \lambda _1( \Sigma ) \varphi  _0=0$ with $ \Vert  \varphi _0  \Vert _{L^2( \Sigma )}=1$. In particular $ \varphi  _0=0$ in $\partial  \Sigma $. Since, as well-known,  the function $H_0(x)=r^ {\alpha_{ \Sigma }} \varphi_0  (\omega )$ --$r= \vert  x \vert  $, $\omega =x/ \vert  x \vert $-- is harmonic,  Theorem \ref{theorem1}  means that any two positive harmonic functions in ${ \mathcal C}_o( \Sigma )$ vanishing on $\partial { \mathcal C}_o( \Sigma )$ are proportional. Note also that Theorem \ref{theorem1} implies that $H_0$ is a positive {\sl minimal} harmonic function in ${ \mathcal C}_o( \Sigma )$.

Section \ref{secttheorem1}  is devoted to a proof of Theorem \ref{theorem1}. See Theorem \ref{theorem2} and an improvement in Theorem \ref{theorem2bis}. It relies in particular on the study of  minimal Martin functions  arising at infinity in ${ \mathcal C}_0( \Sigma )$ and  the study of the convergence in the Martin topology towards such Martin points (for Martin's theory, see \cite{martin}, \cite{naim}, \cite{doob}, or \cite{ancbook}).

 When $ \Sigma $ is  sufficiently regular  Theorem \ref{theorem1} is well-known.  See \cite{kur} for the NTA case. The recent paper  of K. Hirata \cite{hir} establishes the result when $ \Sigma $ is  John. These papers rely on (and provide) Harnack boundary inequalities which do not hold in the general case.

In section \ref{secexample1} we show  that --in contrast with the case where $ \Sigma $ is  John -- another question which might seem at first  to be another formulation  of A.\ Eremenko's question has a negative answer for a general $ \Sigma $, at least in higher dimensions.

\begin{theorem} \label{theoremint2} For $d \geq 4$, there exists a domain $ \Sigma $ such that the Martin boundary of ${ \mathcal C}_o( \Sigma )$ contains a one parameter family of minimal points which are limits of sequences $ \{ P_n \}$  in ${ \mathcal C}_o( \Sigma )$ going to infinity in ${ \mathbb  R}^d$ (and whose all defining sequences go to infinity in ${ \mathbb  R}^d$).
\end{theorem}

 The class of examples provided to prove Theorem \ref{theoremint2}  is strongly related to the construction  by Cranston and McConnell of a bounded domain $ {\mathbb D } $ in ${ \mathbb  R}^3$ with a positive harmonic function  $h$ in $ {\mathbb D } $ such that the lifetime of the $h$-Brownian motion is almost surely infinite \cite{CM}. As shown in \cite{CM}   this cannot happen in a bounded planar domain. There is a corresponding result here given by  the next statement. The only interesting case is  $d=3$.

\begin{theorem}\label{casdleq3co} If $d \leq 3$ and $ \Sigma $ is a domain  in $S_{d-1}$, every sequence  $ \{ P_n \}$  in ${ \mathcal C}_o( \Sigma )$ going to infinity in ${ \mathbb  R}^d$ converges in the Martin topology towards the canonical Martin function $H_0$.  
\end{theorem}

The proof is given in section \ref{sectdim3} and relies on a result of Ba\~nuelos and Davis \cite{bada}. A  similar proof shows that for all $d \geq 2$, every sequence  $ \{ P_n \}$  in ${ \mathcal C}_o( \Sigma )$ going non-tangentially to infinity in ${ \mathcal C}_o( \Sigma )$ converges in the Martin topology to  $H_0$. See Theorem \ref{theorem41} and its proof in section  \ref{sectdim3}. Extensions  based on \cite{ban91} to $d \geq 4$ of Theorem \ref{casdleq3co}  are also mentioned there.

In section \ref{sectionhemisphere}  we supplement Theorem  \ref{theoremint2}  by showing that for $d \geq 4$, $ { \mathcal C}_0(\Sigma) $ may admit plenty of minimal as well as non-minimal Martin  points associated to sequences $ \{ P_n \}$  in ${ \mathcal C}_o( \Sigma )$ going to infinity in ${ \mathbb  R}^d$ (and whose all defining sequences go to infinity in ${ \mathbb  R}^d$). The corresponding Martin functions are of the form $H(x)=r^ \alpha \,  \psi ( { \frac {  x} { \vert  x \vert   }})$ with $-{ \frac {  d-2} {2 }}< \alpha < \alpha _ \Sigma $. This supplementary construction is based on the following fact (see Theorem \ref{theoreminnerhemisphere}) proved also in section \ref{sectionhemisphere} and valid for all $d \geq 2$.
 \begin{proposition} If $ \Sigma $ contains a hemisphere $ \Sigma _+$ of $S_{d-1}$,  then every sequence $ \{ P_n \}$ such that ${ \frac {  P_n} { \vert  P_n \vert   }} \in       \Sigma _+$ and $\vert  P_n \vert   \to  + \infty $ converges towards the canonical Martin point $H_0$.
\end{proposition}

In section \ref{sectjohn}, some implications  on  the  Martin boundary that follow from  regularity conditions  are observed. For example, if $M$ is a  closed John regular subset of $ \Sigma $ (see section \ref{sectjohn}) then every sequence $ \{ P_n \}$ in ${ \mathcal C}_o( M )$ such that $ \vert  P_n \vert   \to   \infty $ converges towards the canonical minimal point $H_0$. This generalizes Hirata's main result in \cite{hir}. 

To prove the above mentioned results, it seems more convenient (and natural) to work with a cylinder model of ${ \mathcal C}_o( \Sigma )$ given by ${ \mathbb  R}\times  \Sigma $  equipped with a suitable elliptic operator  $L$. See section \ref{debut}. We note here that there is a  large literature dealing with Martin's boundaries of product structures (e.g., \cite{molchanov}, \cite{freire}, \cite{murataRIMS},\cite{murata}  or  \cite{GJT}).

Almost all the results (and their proofs) extend to the framework of a cylinder  ${ \mathcal C}_Y( \Sigma ):= { \mathbb  R}\times \Sigma $ where $ \Sigma $ is a relatively compact region in a manifold $M$, the cylinder being equipped with an elliptic operator $L$ in the form $L=(\partial _t \circ \partial_{t}+  \beta \,\partial _t  ) \oplus  { \mathcal L}$ where $ \beta  \in   { \mathbb R} $, $\partial _t$ is differentiation with respect to the first variable and ${ \mathcal L} $ is a second order uniformly elliptic operator in $M$. See section \ref{sectgene}.

\section {Proof of Theorem \ref{theorem1}.}\label{secttheorem1}

{ \em We will  assume once for all that $S_{d-1} \setminus  \Sigma $ is non polar in $S_{d-1}$}. Otherwise, by a standard extension Theorem and the Liouville property, positive harmonic functions in ${ \mathcal C}_0( \Sigma )$ are constant, $ \lambda _1( \Sigma )=0$ and  Theorem \ref{theorem1} is obvious in this case.
Thus, in what follows,  $0$ is a Dirichlet-regular boundary point for ${ \mathcal C}_o( \Sigma )$ -using e.g.\ Wiener's test-, and $ \displaystyle \lim _{x \to  0}h(x)=0$ for the functions $h$ under consideration in Theorem \ref{theorem1}.

\subsection{\normalsize Reduction and some preliminary remarks}\label{debut} As a first step, using a simple change of coordinates  we   reduce the study of positive harmonic functions in the cone 
 $ { \mathcal C}_o( \Sigma ) $  to the study of positive solutions in the cylinder  ${ \mathcal C}_Y( \Sigma )={ \mathbb R} \times  \Sigma $ of some natural  elliptic operator -- the notation ${ \mathcal C}_Y( \Sigma )$ will be used all along the paper--. Recall first the following  classical expression of the Laplacian in polar coordinates: if $f$ is  $C^2$-smooth in the open set  $U$ of ${ \mathbb R} ^d$, $d \geq 2$,
\begin{equation}  \Delta  f(x)= {\frac { \partial ^2 g} {\partial r^2 } }(r,\omega ) +{ \frac { d-1  } {r}}\, { \frac { \partial g  } {\partial r}}(r,\omega )+{ \frac { 1  } {r^2}}\,  (\Delta_S)  _\omega \,g(r,\omega )\end{equation}
 for  $x =r\omega $, $(r,\omega ) \in    \tilde U:= \{ (r,\omega ) \in   { \mathbb R}_+^*  \times S_{d-1}  ;\, r\omega  \in   U\, \}$ and  $g(r,\omega )=f(r\omega) $ for $(r,\omega ) \in    \tilde U$. 
(Recall that $ \Delta _S$ is the Laplace-Beltrami operator on the Riemannian sphere $S_{d-1}$)

Consider then the  change of coordinates:  $\Phi (r,\omega )=  (\log(r), \omega )$, or $ \Phi ^{-1}(u,\omega )=(e^u,\omega )$. Obviously $  \Phi $ defines a $C^ \infty $ diffeomorphism of ${ \mathbb R} _+^*\times  \Sigma $ onto the cylinder ${ \mathcal C}_Y( \Sigma )={ \mathbb R} \times  \Sigma $. By an elementary calculation,
 if $g \in   C^2(  { \mathbb R} _+^*\times  \Sigma )$ and $h=g \circ  \Phi ^{-1}$, we have
\begin{equation}  (\Delta_S)_\omega  g(r,\omega )+ (d-1)\,r\, { \frac { \partial g  } {\partial r}}(r,\omega )+r^2{\frac { \partial ^2g} { \partial r^2} }  = (\Delta_S)_\omega h(u,\omega )+ {\frac {\partial ^2h } {\partial u^2 } }(u,\omega )+(d-2)  {\frac { \partial h} {\partial u } }(u,\omega )\end{equation}

Using these two formulas we are reduced to the study of the Martin boundary of the  cylindrical region ${ \mathcal C}_Y( \Sigma )={ \mathbb R} \times  \Sigma $ of the manifold $X:={ \mathbb R} \times S_{d-1}$ with respect to the second order elliptic operator $ L:=\Delta  _S+(d-2)  {\frac { \partial } {\partial u } }+{\frac { \partial ^2} {\partial u^2 } }$ (where $u$ denotes the ${ \mathbb  R}$-component). 

There is a well known  explicit (and elementary) expression of the heat semi-group $ \{ Q_t \}_{t \geq 0}$ associated to the component  of $L$ acting in ${ \mathbb  R}$, i.e.,  $L_{ \mathbb  R}=\partial ^2_{u} +(d-2)\partial _u$. {\em By heat semi-group associated to $L_{ \mathbb  R}$ } we will mean that for each $ \varphi   \in   C_0^+({ \mathbb  R})$ the function $f: (t,x) \mapsto Q_t( \varphi  )(x)$ is the minimal positive solution of the Cauchy problem : $\partial _tf(t,x)=\partial ^2_xf(t,x) +(d-2)\, \partial _xf(t,x) $ for $t>0$ and $f(0,x)= \varphi  (x)$, $x \in   { \mathbb  R}$.

\begin{lemma} \label{lemme3} For $t>0$ the kernel $  Q_t $ on ${ \mathbb R} $  is given by the convolution with the density \begin{equation}\label{qt} q_t(u)= {\frac {1 } { \sqrt {4 \pi  t}} } \exp(- {\frac {(u+(d-2)t)^2 } {4t } } ),\;\; u \in   { \mathbb  R} 
\end{equation}
\end{lemma}
This means that $Q_t( \varphi  )(u)=\int _{ \mathbb R} q_t(u-v)\,  \varphi  (v)\, dv$ for $   \varphi   \in   C_0({ \mathbb R}  )$.
 The product structure of ${ \mathbb  R}\times  \Sigma $ is exploited through the next standard fact. 

\begin{lemma} \label{lemme4} In $ { \mathbb R}  \times  \Sigma $, the heat semi-group $ \{ H_t \}_{t>0}$ associated to  $ L= {\frac { \partial ^2} {\partial u^2 } } +(d-2) {\frac {\partial  } { \partial u} } +\Delta  _S$ is given by the densities
$$h(t;u,x;v,y)=q_t(u-v)  \pi  _t(x,y),\;\; x,\,y \in    \Sigma ,\;\;u,\,v \in   { \mathbb R},\, t>0. $$
 Here $ \pi  _t $, $t>0$, is the standard heat kernel density in $ \Sigma$ (with respect to $ \Delta  _S$, the usual Riemannian measure $ \sigma  _S$ in $S_{d-1}$ and the Dirichlet boundary condition).
\end{lemma} 

In other words $H_t (\varphi )  (u,x)=\int_{- \infty }^ \infty  \int_ \Sigma \,  h(t;u,x;v,y)\,  \varphi  (v,y)\, dv\, d \sigma_S  (y) $ if $ \varphi   \in   C_0({ \mathbb R} \times  \Sigma ; { \mathbb R} )$ and $(u,x) \in   { \mathcal C}_Y( \Sigma )$.

\subsection {\normalsize Some  inequalities satisfied by Green's function and their consequences}
 The Green's function with pole $( v ,y)$ and with respect to $L$ in ${ \mathbb  R}\times  \Sigma $ (and the measure $\mu  (d(u,x))=du\, d \sigma  (x)$) is the function:
\begin{align}\label{greendef} (u,x) \mapsto G(u,x;v,y)&=\int _0^ \infty h(t;u,x;v,y)\, dt,\;\; (u,x) \in   { \mathbb  R}\times  \Sigma .
   \end{align}

But $h(t;u,x;v,y)=  {\frac {1 } { \sqrt{4 \pi  t} } } \;e^{ -{\frac { (u-v+(d-2)t)^2} {4t } } } \;\pi _t(x,y)$ and for $ \rho >0$ we have
\begin{align}
 e^{ -{\frac { (u+ \rho -v+(d-2)t)^2} {4t } } }  &= e^{ -{\frac { (u -v+(d-2)t)^2} {4t } } }  e^{- {\frac {  \rho } {2t } } ( {\frac { \rho  } {2 } } +(d-2)t+u-v)} \nonumber \\
{\hspace{2.7truecm} }  & \geq  e^{ -(d-2)\rho /2} \,e^{ -{\frac { (u -v+(d-2)t)^2} {4t } } } 
\end{align}
if $v \geq u+{ \frac {  \rho   } {2}}$.

On the other hand, if $ \rho >0$, and $v \leq u+{ \frac {  \rho   } {2}}$
\begin{eqnarray}
 e^{ -{\frac { (u+ \rho -v+(d-2)t)^2} {4t } } }  = e^{ -{\frac { (u -v+(d-2)t)^2} {4t } } }  e^{ {\frac {  \rho } {2t } } ( -{\frac { \rho  } {2 } } -(d-2)t+v-u)} \leq e^{ -(d-2)\rho/2 } \,e^{ -{\frac { (u -v+(d-2)t)^2} {4t } } } \nonumber \\ 
\end{eqnarray}
  Upon integrating with respect to $t$ these inequalities, we obtain the next proposition.

\begin{proposition} \label{proposition1} The Green's function $G$ for $L$ in ${ \mathbb R} \times  \Sigma $ satisfies the following relations: 
$$G(u,x;v,y)  \leq e^{(d-2) \rho/2 } \,G(u+ \rho ,x;v,y), \; \; {\rm \ for \ }u,\,v \in   { \mathbb  R},\;\, {\rm and \ } x,\,y \in    \Sigma ,$$ if $ \rho >0$ and $v\geq  {\frac { \rho  } { 2} } +u$. And $G(u,x;v,y)  \geq e^{(d-2) \rho/2 } \,G(u+ \rho ,x;v,y)$ if $ \rho >0$ and $v\leq  {\frac { \rho  } { 2} } +u$.

Moreover, we have the following symmetry identities :  for $\; \; x,y \in    \Sigma, u, v,v_0, v_1 \in   { \mathbb  R}$
\begin{align}\label{symm}G(v_0-u,x;v_0-v,y)=e^{(d-2)(u-v)}\,G(v_1+u,x;v_1+v,y),  
\end{align}

\end{proposition}

Note that (\ref{symm}) follows immediately from (\ref{qt}), lemma  \ref{lemme4} and (\ref{greendef}). Of course these identities can be understood in terms  of Kelvin's transformation  if one returns to the cone ${ \mathcal C}_o( \Sigma )$ equipped with the usual Laplacian.

The above leads to the following properties of the $L$-Martin functions associated to $L$ in the cylinder ${ \mathcal C}_Y( \Sigma )$. We choose and fix once for all a reference point $x_0 \in    \Sigma $ and take $(0,x_0)$ as the normalization point for Martin's functions in ${ \mathcal C}_Y( \Sigma )$. Occasionally, we use the standard notations ${\bm  \Delta  }$ (resp.\ ${\bm  \Delta  }_1$) to denote the Martin boundary (resp.\ the minimal Martin boundary) of $({ \mathcal C}_Y( \Sigma ), L)$,  and $ \widehat { \mathcal C}_Y( \Sigma )$ its Martin compactification (ref. \cite{martin}, \cite{naim}, \cite{ancbook}).

\begin{proposition} \label{proposition2}
If $K$ is an $L$-Martin function in ${ \mathcal C}_Y( \Sigma )={ \mathbb  R}\times  \Sigma $ defined by a sequence $(v_j,y_j)$ with  $v_j \to  + \infty $, $y_j \in    \Sigma $, (i.e., $K(u,x)=\lim_{j \to   \infty } K_{(v_j,y_j)}(u,x)$ where $K_{(v_j,y_j)}$ is the Martin kernel $K_{(v_j,y_j)}:=G(.,.;v_j,y_j)/G(0,x_0;v_j,y_j)$) then  
\begin{equation}\label{eqproposition2}K(u+ \rho ,x) \geq e^ {-(d-2)\rho/2} K(u,x),\;\;\; (u,x) \in   { \mathbb R} \times  \Sigma \end{equation}
 for every $ \rho >0$. In particular if $K$ is minimal there exists  $ \alpha  \geq   -{\frac { d-2} { 2} } $ such that  $K(u,x)=e^{ \alpha u} K(0,x)$,   $(u,x) \in   { \mathbb  R}\times    \Sigma   $, and the function $ s(x)=K(0,x)$ is a proper function in $ \Sigma $, namely $ \Delta_S(  s)+ \lambda s=0$ for $ \lambda = \alpha ^2+(d-2) \alpha $. So $ \lambda  \geq -{ \frac { (d-2)^2  } {4}}$, $ \alpha ={ \frac {2-d+\sqrt{(d-2)^2+ 4\lambda }   } {2}} $ and $s$ is a minimal positive $( \Delta  _S+ \lambda I)$--harmonic function in $ \Sigma $.

\end{proposition}

The first statement follows from Proposition \ref{proposition1} and the definition of Martin functions. If $K$ is minimal, then $(u,x) \mapsto K(u+ \rho ,x)$ is also $L$-minimal; by (\ref{eqproposition2}), it follows that $K(u+ \rho ,x)=c( \rho )K(u,x)$ for some function $c \in   C({ \mathbb  R};{ \mathbb  R}_+^*)$ and all $( \rho ,u,x) \in   { \mathbb  R}\times { \mathbb  R}\times  \Sigma $ (assuming first $ \rho >0$). We have  $c( \rho + \rho ')=c( \rho )\, c( \rho ')$ for $ \rho ,\,  \rho ' \in   { \mathbb  R}$ and so  $c( \rho )=e^{ \alpha  \rho }$ for some $ \alpha  \geq -{ \frac { d-2  } {2}}$. The other claims are then immediate. $\square$

Note that $ \alpha  \in   [ \alpha _0,  \alpha _ {\rm max} ]$ where $ \alpha _{\max}:= \alpha _ \Sigma ={ \frac {2-d+\sqrt{(d-2)^2+ 4\lambda_1 }   } {2}} $ and $ \alpha _0=-{ \frac {  d-2} {2 }}$. We set $ \alpha _{\rm min}={ \frac {2-d-\sqrt{(d-2)^2+ 4\lambda_1 }   } {2}} $ .

\begin{remark}\label{remark2} Similarly if $K$ is an $L$-Martin function related to a sequence $(v_j,y_j)$, with $v_j \to  - \infty $, $y_j \in    \Sigma $, then $K(u+ \rho ,x) \leq e^ {-(d-2)\rho/2} K(u,x)$ when $ (u,x) \in   { \mathbb R} \times  \Sigma$ and $ \rho  \geq 0$. If $K$ is minimal w.r.\ to $L$, then $K(u,x)=e^{ \alpha u}s(x)$, for some $ \alpha   \in   [ \alpha _{\min},  \alpha _0] $ and some $ \lambda $-proper function $s$ in $ \Sigma $, i.e., $ \Delta  _S(s)+ \lambda s=0$ in $ \Sigma $. Here $ \lambda = \alpha ^2+(d-2) \alpha  \geq -{ \frac { (d-2)^2  } {4}}$,  $ \alpha ={ \frac {2-d-\sqrt{(d-2)^2+ 4\lambda }   } {2}} $ and  $s$ is $ (\Delta  _S+ \lambda I)$-minimal in $ \Sigma $.\end{remark}

Observe   that if ${\bm  \Delta  }^+$ (resp. ${\bm  \Delta  }^-$) is the set  of the Martin functions arising from a sequence $(v_j,y_j)$ with $\lim y_j=+ \infty $ (resp.  $\lim y_j=- \infty $) the identities (\ref{symm})  lead to a natural bijection $K \mapsto   \tilde K$ from ${\bm  \Delta  }^+$ onto ${\bm  \Delta  }^-$, where $\tilde K(u,x):=e^{-(d-2)u}K(-u,x)$.

We  note two other  observations which complement Proposition \ref{proposition2} and Remark \ref{remark2}.

\begin{proposition} \label{proposition4} If $ \zeta   \in    {\bm  \Delta  } $ is a Martin boundary point for $({ \mathcal C}_Y( \Sigma ),L)$ which  as a point of the Martin compactification is in the  closure of both  $ \{ (u,x) \in   { \mathcal C}_Y( \Sigma )\,;\, u \geq  \rho \,  \}$ and $ \{ (u,x) \in   { \mathcal C}_Y( \Sigma )\,;\, u  \leq   \rho' \,  \}$ for some $ \rho,\,  \rho '  \in   { \mathbb  R}$, $ \rho '< \rho $, then $K_ \zeta  $ is in the form $K_  \zeta  (u,x)=e^{-{ \frac { d-2  } {2}}u}\,f(x)$ where $f$ is a positive solution of $ \Delta  _Sf-{ \frac { (d-2)^2  } {4}} f=0$ in $ \Sigma $.

The same conclusion holds for $ \zeta   \in    {\bm  \Delta  }_1$, if $ \zeta  $ is in the closure (for the Martin topology) of a set $ { \mathcal C}_Y^R=\{ (v,y)\,;\,  \vert  v \vert   \leq R,\, y \in    \Sigma \,  \}$, $R>0$, and if $K_ \zeta  $ is of the form $K_ \zeta  (u,x)=e^{ \alpha u} f(x)$. Moreover $f$ is a positive  minimal solution of $ \Delta  _Sf-{ \frac { (d-2)^2  } {4}} f=0$ in $ \Sigma $. 
\end{proposition}

In the first case, by the above estimates of Green's function we must have $K_  \zeta  (u,x)=e^{-{ \frac { d-2  } {2}}u}\,f(x)$ for $ \rho '<u< \rho$, $x \in    \Sigma $, where $f(x)=e^{{ \frac { d-2  } {2}} \rho '}\,K_ \zeta  ( \rho ',x)=e^{{ \frac { d-2  } {2}} \rho }\,K_ \zeta  ( \rho ,x)$. Then $f$ is  necessarily as in the statement and by the Green's function estimates, we see that $K_  \zeta    (u,x)   \geq  e^{-{ \frac { d-2  } {2}}u}\,f(x)$ for all $u> \rho '$. 

But  a nonnegative solution $w$ of $L(w)=0$ in a domain $\Omega  \subset   { \mathcal C_Y}( \Sigma )$ that vanishes at some point vanishes everywhere (by Harnack inequalities). Thus $K_  \zeta    (u,x)   =  e^{-{ \frac { d-2  } {2}}u}\,f(x)$ when $u > \rho '$. A similar argument extends  this equality to $u< \rho $. 

In the second case, if $ \zeta  =:\lim (v_j,y_j)$, we must also have $ \zeta  =\lim (v_j+s,y_j)$ for every $s \in   { \mathbb  R}$. So the result follows from the first part of the proposition (the minimality of $f$ being necessary for the minimality of $K_ \zeta  $). $\square$


\subsection{\normalsize Positive $L$-harmonic functions vanishing on ${ \mathbb  R}\times  \partial \Sigma $ (a)} 
Let  $ f$ denote a nonnegative $L$-harmonic function  in ${ \mathcal C}_Y( \Sigma ):={ \mathbb R} \times  \Sigma $ such that  $ f=0$ on ${ \mathbb R} \times \partial  \Sigma$ in the  weak sense. 

By the Martin boundary  theory there is a unique integral representation of $f$ in the form
\begin{equation} f(u,x)=\int K_ \zeta  (u,x)\, d\mu  ( \zeta  ),\;\;\;(u,x) \in    { \mathbb R}  \times  \Sigma 
\end{equation} where $\mu  $ is a finite positive Borel measure on the minimal Martin boundary ${\bm  \Delta  }_1$ of $({ \mathcal C}_Y( \Sigma ),L)$ and where $K$  is the $L$-Martin kernel associated to $L$ and the reference point $(0,x_0) \in    \Sigma $.

Denote ${\bm  \Delta  }_1^ \infty $ the trace on ${\bm  \Delta  }_1$ of the intersection of the closures --w.\,r.\ to Martin's topology-- of the sets $ X_N:= \{ (- \infty ,-N] \cup [N, +\infty ) \}\times  \Sigma $, $N \geq 1$. In other words, a point $  \zeta   \in   {\bm  \Delta  }_1$ is in ${\bm  \Delta  }_1^ \infty $ if and only if there is a sequence of points $ \xi  _j=(u_j,x_j) \in   { \mathbb R} \times  \Sigma $ with $ \vert  u_j  \vert  \to  + \infty$  converging to $ \zeta  $. In particular $K_ \zeta  $ is in the form given by either Proposition \ref{proposition1} or remark \nolinebreak \ref{remark2}.

\begin{lemma} The measure $\mu_f  $ is supported by ${\bm  \Delta  }_1^ {\infty} $.
\end{lemma}
Proof. Let $N \geq 1$. In ${ \mathcal C}_Y( \Sigma )$,  the function $f$ is equal to $R_f^{X_N}$, its  r\'{e}duite (cf.\ \cite{brelot}, \cite{CC}, \cite{doob})  on $X_N$    (w.r.\ to $L$). This follows from the assumption $f=0$ in $\partial  \Sigma \times { \mathbb  R}$ and from (a standard extended form of) the maximum principle. So by  the Martin boundary theory, the measure $\mu  _f$ is supported by the set of the points $  \zeta     \in   {\bm  \Delta  }_1$ such that $X_N$ is not minimally thin at $ \zeta  $. Such a point $ \zeta  $ is necessarily in the closure of $X_N$ in $ \widehat { \mathcal C}_Y( \Sigma )$, the $L$-Martin compactification of ${ \mathcal C}_Y( \Sigma )$. Whence the result. $\square$

Next we will observe  a simple condition for $\mu  _f$ to be concentrated on 
\begin{equation}{\bm  \Delta  }_1^{+ \infty } := {\bm  \Delta  }_1 \cap \{  \zeta   \in    \widehat { \mathcal C}_Y( \Sigma )\,;\,  \zeta  =\lim_j (u_j,x_j) {\rm \ with\ } \;\,x_j  \in    \Sigma \, {\rm \ and \ } u_j  \to   + \infty\;   \}. \end{equation} 

Define similarly $\Delta  _1^{- \infty }$ using sequences $ \{ (u_j,x_j) \}_{j \geq 1}$ with $\lim u_j=- \infty $.

\begin{lemma} \label{lemma6}  If $f(-t,x_1)=o(e^{ {\frac { d-2} { 2} } t})$ as $t \to  + \infty $ for some (or all) $x_1 \in    \Sigma $, then the measure $\mu_f  $ is supported by ${\bm  \Delta  }_1^ {+\infty} $. 
\end{lemma} 

This is because for $u \leq 0$ and $ \zeta   \in  {\bm  \Delta  }_1^{- \infty }$ we have  $K_ \zeta  (u,x_0) \geq e^{- {\frac {d-2 } { 2} }u} \,K_ \zeta  (0,x_0)=e^{- {\frac {d-2 } { 2} }u}$ (see Remark \ref{remark2}) and 
$$f(u,x_0)=\int K_ \zeta  (u,x_0)\, d\mu  ( \zeta  )  \geq  e^{- {\frac {d-2 } { 2} }u} \mu  ({\bm  \Delta  }_1 ^\infty \setminus {\bm  \Delta  }_1^{+ \infty } )=e^{- {\frac {d-2 } { 2} }u} \mu  (  {\bm  \Delta  }_1^{- \infty } ).$$

Thus $\mu  ({\bm  \Delta  }_1^{- \infty })=0$. Note the special case where $f(t,x_0) =O(1) $ as $t \to  - \infty $. $\square$


\subsection{ \normalsize End of proof of Theorem \ref{theorem1}} 

In this subsection it is assumed that {\em moreover  $f(u_n,x_0)=O(1)$ for a real sequence $u_n \to  - \infty $}, 
that is   $ \displaystyle \liminf_{u \to  - \infty } f(u,x)< \infty $ for each $x\in    \Sigma $. We will show that  $f$ is unique up to a multiplication by a constant and will thus prove Theorem \ref{theorem1}. Examining the proof, it will be seen later that this assumption on $ \displaystyle \liminf_{u \to  - \infty } f(u,x_0)$ can be notably relaxed.

For any given real $ a _0$ let us  denote ${\bm  \Delta  }_1^ {+\infty} ( a _0)$   the set of all points  $  \zeta   \in   {\bm  \Delta  }_1^ {+\infty} $ such that $K_ \zeta  $ is in the form : $K_ \zeta (u,x) =  \varphi(x)  \, e^{  \alpha  u}$, $(u,x) \in   { \mathcal C}_Y( \Sigma )$ for some real $ \alpha  \geq  a _0$ and some   minimal positive solution of $ \Delta  _S \varphi  + \lambda  \varphi  =0$ in $ \Sigma $ where $  \lambda = \alpha ^2+(d-2) \alpha $. We set $ \alpha = \alpha ( \zeta  )$ and $ \lambda = \lambda ( \zeta  )$.

Repeating the argument in lemma \ref{lemma6} we first note the following.

\begin{lemma} The measure $\mu  _f$ is concentrated  in ${\bm  \Delta  }_1^ {+\infty} (0)$.
\end{lemma} 

This is again immediate since given $ \varepsilon >0$ we have   for $u \leq 0$,    
$$f(u,x_0) \geq \int _{{\bm  \Delta  }_1^{+ \infty} \setminus {\bm  \Delta  }_1^ {+\infty} (- \varepsilon )} K_  \zeta    (u,x_0)\, d\mu  _f( \zeta  ) \geq e^{- \varepsilon u}\mu  _f({\bm  \Delta  }_1^ \infty \setminus {\bm  \Delta  }_1^ \infty (- \varepsilon  ))\;\, $$
and hence $\mu  _f({\bm  \Delta  }_1^{+\infty} \setminus {\bm  \Delta  }_1^{+ \infty }(- \varepsilon ))=0$ since $ \displaystyle \liminf _{u \to  - \infty }f(u,x_0)< \infty $. $\square$
 
\medskip  
Recall that $ \lambda _1=\lambda _{\max}$ is the greatest  $ \lambda $ for which the equation $ \Delta  _S \varphi  + \lambda  \varphi  =0$ admits a positive solution (or even a positive supersolution) in $ \Sigma $.

\begin{lemma}\label{lemmamuf(A)vanish}  Set $A= \{  \zeta   \in   {\bm  \Delta  }_1^ {+\infty} (0)\,;\, \lambda ( \zeta  )< \lambda _{\max}( \Sigma )\,  \}$. Then $\mu  _f(A)=0$.
\end{lemma} 

Assume that $\mu  _f(A)>0$. Then for some $ \lambda '_1$, $0< \lambda '_1< \lambda _1$, the set $A'$ of all $ \zeta   \in   {\bm  \Delta  }_1^{+ \infty }(0)$ for which  $ 0 \leq \lambda ( \zeta  ) \leq  \lambda '_1$ has strictly positive measure: $\mu  _f(A')>0$. Note that $A'$ is the set of the minimal points $ \zeta   \in   {\bm  \Delta  }_1$ such that for some $ \lambda $, $0 \leq  \lambda  \leq  \lambda '_1$, 
 \begin{equation}  K_{ \zeta  }(u,x)=  e^{ \alpha u}\, k_ \zeta    ^ \lambda (x)  \end{equation}
where $\alpha = {\frac { -(d-2)+\sqrt{(d-2)^2+4 \lambda }} { 2} }  $  and $ k_ \zeta  ^ \lambda  $ is a $\Delta  _S+ \lambda I$--minimal  function  in $  \Sigma $. 

The function 
\begin{equation}\label{IRphi} \varphi  (x)=\int _{A'} k_  \zeta    ^  \lambda (x)\, d\mu  _f(  \zeta  ), \;\; x \in    \Sigma ,
\end{equation}

is a (strictly) positive $ \Delta  _S$-superharmonic function in $ \Sigma $ which satisfies the following:

\ \  (i) $ \varphi  $ vanishes weakly  on the boundary of $ \Sigma $ 
(note that $ \varphi  (x) \leq f(u,x)$ for $x \in    \Sigma $, $u \geq 0$ and that by assumption $f=0$  on ${ \mathbb R}  \times \partial  \Sigma $ in the weak sense), 

\ \ (ii) the positive measure $ - \Delta_S  (\varphi)  $ admits  the density  $ \psi (x)  =\int_{A'}  \lambda k_  \zeta    ^  \lambda (x)\, d\mu  _f(x) $ with respect to $  \sigma  _{S}$. 

By (i) and a well-known form of the maximum principle $ \varphi  $ is a potential in $ \Sigma $ with respect to the spherical Laplacian $ \Delta  _S$ in $\Sigma $, i.e.\ its greatest subharmonic minorant $h$ in $ \Sigma $  is zero . Indeed there exists  a positive $ \Delta  _S$-superharmonic function $s$ going to infinity at each $ \zeta   \in   \partial  \Sigma $ where $  \varphi  $ has a non-zero upper limit since the set of these points $ \zeta  $ is polar.  So by a standard form of the  maximum principle $h- \varepsilon s \leq 0$ for each $ \varepsilon >0$. 

It follows that $  \varphi  $ is a Green's potential in $ \Sigma $ and so using first (ii) and then (\ref{IRphi}) and (ii) again, we get
\begin{align} \label{eqGsigmapsi}
   G_ \Sigma ( \psi )&= \varphi   \geq   {\frac { 1} { \lambda '_1 } }  \psi  \end{align} 
in $ \Sigma $, where $G_ \Sigma$ denotes the Green's function in $ \Sigma $ w.r.\  to $ \Delta  _S$. 

The function $ \psi  $ is in $   L^2( \Sigma  )$ (it is bounded in $ \Sigma $ -in fact $ \psi $ is also an element of $H_0^1( \Sigma )$-). Thus from $ \psi \ne 0$ and $ \Vert G_ \Sigma ( \psi ) \Vert _2 \geq  {\frac { 1} {  \lambda' _1} }  \Vert  \psi  \Vert _2$   we infer that $ {\frac {1 } {  \lambda _1} }  =  \Vert G_ \Sigma  \Vert _{L(L^2,L^2)} \geq  {\frac {1 } { \lambda' _1 } } $ which is absurd. $\square$

On the other hand, we know that  every positive eigenfunction $ \varphi   $ with respect to $( \Delta_S  , \Sigma )$ and the eigenvalue $  \lambda _1= \lambda _1( \Sigma )$ is proportional to $ \varphi  _0$ (see e.g.\ \cite{ancbook}). Moreover $ \varphi  _0$ vanishes at every Dirichlet-regular boundary point $\xi   \in   \partial _S \Sigma $ and so vanishes quasi-everywhere on $\partial _S \Sigma $.

Thus we have proved  the following which contains Theorem \ref{theorem1}. See also Theorem \ref{theorem2bis}.
\begin{theorem}\label{theorem2} There is  a unique positive $L$-harmonic function $F_+$ on ${ \mathcal C}_Y( \Sigma )={ \mathbb R} \times  \Sigma $ vanishing (in the weak sense) on the boundary ${ \mathbb R} \times \partial  \Sigma $ and such that $ \displaystyle \liminf_{u \to  - \infty}  F_+(u,x_0)< +\infty $ and $F_+(0,x_0)=1$. This function is a minimal Martin function  and $F_+(u,x)=e^{ \alpha_ {\max}  u}  {\frac { \varphi  _0(x)} {\varphi  _0(x_0) } }  $ if $ \alpha_ {\max}  = {\frac {2-d+\sqrt {(d-2)^2+4 \lambda _1 }} { 2} } $ (i.e.\,$ \alpha _{\max}= \alpha _ \Sigma $).

\end{theorem}

  In what follows we will keep the notation $F_+$ for this ``canonical" minimal function and let $F_-(u,x)=e^{ \alpha_{\min} u}  {\frac {\varphi  _0(x)} {\varphi  _0(x_0) } }$, where $ \alpha_{\min} ={\frac {2-d-\sqrt {(d-2)^2+4 \lambda _1 }} { 2} } $ for the similar minimal function related to the end $``u \to  - \infty "$ of the cylinder. We set $ \alpha _0=-{ \frac { d-2  } {2}}$.

In the course of the proof of Theorem \ref{theorem2} we have  also essentially shown the following facts. Distinguish three classes of minimal Martin functions $K$ w.r.\ to  $({ \mathcal C}_Y( \Sigma ),L)$: \linebreak (i) the functions in the form $K_  \zeta    (u,x)=k(x)\, e^{ \alpha u}$  with  $   0<\vert  \alpha -\alpha _0 \vert   \leq {\frac {\sqrt {(d-2)^2+4 \lambda _1 }} { 2} } $, (ii) those in  the form $K_  \zeta    (u,x)=k(x)\, e^{ \alpha _0 u}$,   (iii)  the class of all other minimal functions.  

\begin{proposition} \label{classif}  If $K_ \zeta  $ is in  the third class there is a   unique $v_ \zeta   \in   { \mathbb R} $ such that $\lim v_j=v_ \zeta  $ for any sequence $ \{ (v_j,y_j) \}$ in ${ \mathcal C}_Y( \Sigma )$ converging to $ \zeta  $ and for such a sequence  $\lim_j(w_j,y_j) =  \zeta  $ whenever $w_j \to  v_ \zeta  $ in ${ \mathbb R} $. If  $K_ \zeta  $ is in the first class, every corresponding sequence $ \{ (v_j,y_j) \}_{j \geq 1}$ is such that $\lim v_j=+ \infty $ or $\lim v_j=- \infty $  depending whether $ \alpha > \alpha _0$ or $ \alpha < \alpha _0$. Finally if $K_ \zeta  $ is in the  class (ii) there are sequences $(v,y_j)$ (with a fixed first coordinate) converging to $ \zeta  $ and for any such sequence $\lim_j(v_j,y_j) =  \zeta  $ for every real bounded sequence $ \{ v_j \}$; moreover there are sequences $ \{ v_j \}$ such that $\lim v_j=+ \infty $ and $\lim_j(v_j,y_j) =  \zeta  $.
\end{proposition}

A minimal function in the  class (iii)  will be said to be of the finite type. 

Of course if $ \Sigma $ is smooth, the first class reduces to $ \{ F_+,F_- \}$ and the second class is empty. We shall see later  that  there may exist minimal as well as non minimal Martin points $ \zeta  $  in the form $K_ \zeta  (u,x)=e^{ \alpha u}k(x)$, for all $ \alpha $ such that $ \vert   \alpha - \alpha _0 \vert  < {\frac {\sqrt {(d-2)^2+4 \lambda _1 }} { 2} } $. See \ref{sectionhemisphere}.

{\bf Proof.} To establish the last claim let $  \zeta  $ be  in the second class. If $v \in   { \mathbb R} $  and if $(v_j,y_j)  \to   \zeta  $ then  $(2v-v_j,y_j) \to   \zeta  $ (by the identities (\ref{symm}) with $v_0=v_1=0$). But a minimal Martin point has a neighborhood basis $ \{U_j \}$ in $ \widehat { \mathcal C}_Y( \Sigma )$ with $U_j \cap    { \mathcal C}_Y (\Sigma )$ connected (by the general theory, see e.g.\ \cite{naim} \ p.\ 223) and so  we can find points $z_j  \in    \Sigma $ with $(v,z_j) \to   \zeta  $. If $ \{ v_j \}$ is bounded it follows at once from the (local) Harnack inequalities and the translation invariance with respect to the first coordinate that $(v_j,z_j) \to   \zeta  $. It is then obvious that if $v_j \to  + \infty $ sufficiently slowly $(v_j,y_j) \to   \zeta $. 

It also follows immediately from translation invariance that  if a sequence $ \{ (v_j,y_j) \}_{j \geq 1}$ converges to a point $ \zeta  \in    \Delta     $  then $\lim_j (w_j,y_j) =  \zeta  $ for $ \{ w_j \}$ such that $ \vert  v_j-w_j \vert    \to  0$.  $\square$

\subsection{\normalsize Positive $L$-harmonic functions vanishing on ${ \mathbb  R}\times  \partial \Sigma $ (b)} 

The proof of Lemma \ref{lemmamuf(A)vanish} can be extended so as to use a  much weaker assumption on the behavior of $f(u,x_0)$ for $u \to  - \infty $. This leads to a description  of the positive $L$-solution in $ { \mathbb R} \times \Sigma $ vanishing on ${ \mathbb R} \times  \partial \Sigma $ which also improves Theorem \ref{theorem2}.

\begin{theorem}\label{theorem2bis} If $f$ is $L$-harmonic in ${ \mathcal C}_Y( \Sigma )$ and vanishes on ${ \mathbb  R}\times \partial  \Sigma $, then $f$ is a linear combination of $F_+$ and $F_-$. Thus if  $ \displaystyle \liminf_{u \to  - \infty } e^{ \alpha _{\rm min} u} f(u,x_0)=0$ then $f$ is proportional to $F_+$.
\end{theorem}

Using the Martin disintegration of $F$, we may write $F=aF_++bF_-+F_1+F_2$ with 
\begin{align}
  F_1(u,x)=\int_{A} e^{ \alpha u}\,k_ \zeta  (x)\, d\mu  ( \zeta  ) &\;,\; F_2(u,x)=\int_{B} e^{ \alpha u}\,k_ \zeta  (x)\, d\nu  ( \zeta  )\nonumber \end{align}
 where $A= \{  \zeta   \in   {\bm  \Delta  }_1\,;\,  \exists   \alpha _ \zeta  $,  $\alpha _0 \leq  \alpha_ \zeta   < \alpha _{\rm max} $ and $k_ \zeta   \in   C_+( \Sigma )$ s.t.\ $K_ \zeta (u,x) \equiv e^{ \alpha _ \zeta  u}\, k_ \zeta  (x) \,  \,  \}$, $B=\{  \zeta   \in   {\bm  \Delta  }_1\,;\,  \exists   \alpha _ \zeta  $,  $\, \alpha _{\rm min}<\alpha_ \zeta    <  \alpha _0\, $ and $k_ \zeta   \in   C_+( \Sigma )$ with $K_ \zeta (u,x) \equiv e^{ \alpha _ \zeta  u}\, k_ \zeta  (x) \,  \,  \}$ and where $\mu  $ and $\nu $ are finite Borel measures supported by $A$ and $B$ respectively. 

We claim that $\mu(A) =0$. If not there exists $ \alpha '_1 \in   ( \alpha _0, \alpha _{\rm max})$ such that $\mu  (A')>0$ if $A'= \{  \zeta   \in    A\,;\,   \alpha _0 \leq  \alpha_ \zeta    \leq  \alpha'_1  \}$ and  repeating the argument in Lemma  \ref{lemmamuf(A)vanish} we may conclude 
using now  $ \varphi  (x):=\int _{A'} k_ \zeta  (x)\, d\mu  ( \zeta  )$ and  the potential theory w.\ r.\ to the operator $ \Delta  _S-{ \frac {  (d-2)^2} { 4}}I= \Delta  _S- \lambda _0I$ in $S_{d-1}$ (in particular the related Green kernel $G_ \Sigma ^{  \lambda  _0}$). Note that (\ref{eqGsigmapsi}) becomes : $G_ \Sigma ^{  \lambda  _0}( \psi )= \varphi   \geq   {\frac { 1} { (\lambda '_1 - \lambda _0)} }  \psi   $ where $ \lambda '_1$ is the eigenvalue corresponding to $ \alpha '$.

In the same way (or using the observation after remark \ref{remark2}) it is shown that $F_2=0$. $\square$

\section{Nontangential convergence  to $F_+$ or $H_0$}\label{sectnontangential}

The next statement is about how Martin's topology relates to the canonical minimal $F_+$ in the cylinder ${ \mathcal C}_Y( \Sigma )$ (or to the minimal $H_0$ in the cone $ { \mathcal C}_0( \Sigma )$). It says   that  nontangential convergence of the current point  $(u,x) \in   { \mathcal C}_Y( \Sigma )$ (resp.\  $x=r\omega  \in   { \mathcal C}_o( \Sigma )$) to the end ``$u = + \infty $" (or ``$ r=  + \infty $") implies its convergence to the canonical Martin  point at infinity.


\begin {theorem} \label{theorem41} For every sequence $ \zeta  _j:=(u_j,x_j)$  in ${ \mathbb R} \times  \Sigma $ such that  $u_j \to   +\infty $ and $ \{ x_j \}$ is relatively compact in $ \Sigma $, it holds  that  $K_{ \zeta  _j}(u,x) \to  e^{ \alpha _ {\max} u} \,\varphi  _0(x)/ \varphi  _0(x_0)$ (i.e. $ \{  \zeta  _j \}$ converges to the Martin function $F_+$). In fact, the following Harnack boundary inequalities hold
\begin{align}\label{eq41}  C^{-1}\,G(u,x_0; v,x_0)\, G(v,x_0; w,x_0) \leq G(u,x_0; w,x_0) \leq C\,G(u,x_0; v,x_0)\, G(v,x_0; w,x_0)  
\end{align} for $u,v,w \in   { \mathbb  R}$, $u+1 \leq v$, $v+1 \leq w$ and some constant $C=C(d, \Sigma, x_0 ) \geq 1$.  The  inequalities obtained by replacing $G$ by its  transposed kernel in (\ref{eq41}) also hold.
\end{theorem}  

A proof will be given at the end of  section \ref{sectdim3}. We note here that inequalities (\ref{eq41}) imply by themselves that for $t \to   +\infty $ the point  $(t,x_0)$ converges to a minimal point in the Martin boundary (see \cite{anc0} Th\'{e}or{\`e}me 2 or \cite{anc1} p.\ 516).

\section{\hspace{-2truemm} Martin boundary  and   subsets of $ \Sigma $} \label{sectjohn}
In this section we  collect   some properties of  the Martin boundary of ${ \mathcal C}_Y( \Sigma )$ resulting from   regularity conditions on a subset $M$ of $ \Sigma $. The results will not be used before section \ref{example1description}.

\subsection{\normalsize John regular subsets}\label{sectjohnreg}
Let $M$ be  a closed subset of $ \Sigma $. For $  \eta   >0$, let $M _  \eta := \{ x \in    \Sigma \,;\, d_a(x,M) \leq   \eta   \,  \}$ where  $d_a(x,y):=\inf  \{ {\rm diam} (C)\,;\, C \subset    \Sigma $ connected, $x,\, y \in    \Sigma \,  \}$ - $S_{d-1}$ is equipped with the usual metric in ${ \mathbb R} ^d$-. For $0<c_0 \leq 1$ we say that $M$ is $c_0$-John in $ \Sigma $ if there are points $A_1$, $\dots$, $A_N$ in $ \Sigma $  such that  (i) $d(A_j, \partial  \Sigma ) \leq c_0^{-1} d(A_k, \partial  \Sigma )$ for $1 \leq j,\, k \leq N$,   (ii) for  $  \eta :=c_0\,\max  \{ d(A_j,  \partial  \Sigma )\,;\, 1 \leq j \leq N\,\}$,   each  $a \in M_ \eta  $  can be connected to one   $A_j$ by a $c_0$-John arc  in $ \Sigma $  (see \cite{anc4} D\'{e}finition 1.1 and Th\'{e}or\`{e}me 5.3). $ \{ A_j \}_{1 \leq j \leq N}$ is   then called  a $c_0$-admissible set of poles for $M$ (note that $N$ can always be chosen smaller than  a constant $N_0(d,c_0)$). 

The next statement generalizes  Hirata's main result in \cite{hir}. We rely on Theorem \ref{theorem41} and  a boundary Harnack principle  given in \cite{anc4} (\cite{AHL} for $N \leq 2$). Note  that  this statement  may be  easily reduced to the $N=1$ case.

\begin{theorem} \label{presdubord} Let $M $ be a closed and $c_0$-John subset of $ \Sigma $ with poles $A_j$, $1 \leq j \leq N$. Then  $lim_{u \to  + \infty }(v,y)=F_+$  in the Martin topology, uniformly with respect to $y \in   M$.
\end{theorem}

{\bf Proof of Theorem \ref{presdubord} } Denote $K$  the Martin kernel in $ \widehat { \mathcal C}_Y( \Sigma )$ with respect to the normalization point $(0,x_0)$.
 Applying  Th\'{e}or\`{e}me 5.3 and Remarque 5.4  in \cite{anc4}   to $M\times [v-1, v+1]$ as a subset of ${ \mathcal C}_Y( \Sigma )$ (or rather -so as to deal with the classical Laplacian- to the corresponding situation in ${ \mathcal C}_0( \Sigma )$) we obtain a constant  $C=C(d,c_0) \geq 1$ such that 
\begin{align}\label{HJB} K_{(v,y)}(u,x) \leq C\, \sum _{j=1}^N \;{ \frac {K_{(v,y)}(v,A_j)  } {K_{(v,A_j)} (v,A'_j)}}\;K_{(v,A_j)}(u,x) \end{align} 
whenever $(v,y) \in   { \mathbb  R}\times M$, $ \vert  v \vert   \geq 1$,  and $(u,x) \in   { \mathcal C}_Y( \Sigma )$ satisfies $ \vert  u-v \vert   \geq 1$ (or $d_a(x,M) \geq  c_0 $). Here $A'_j$ is  arbitrarily chosen in $\partial B(A_j,{\frac  { c_0} {100}} \,{\rm dist}(A_j;S_{d-1}\setminus  \Sigma ))$ and we restrict to $y$ such that $ \vert  y-A_j \vert   \geq 2\,  \vert  A_j-A'_j \vert  $ for all $j$.

By Harnack inequalities  \begin{align}   K_{(v,y)} (u,x) \geq c\;{ \frac {K_{(v,y)}(v,A_j)  } {K_{(v,A_j)} (v,A'_j)}} K_{(v,A_j)}(u,x)
 &\nonumber
\end{align}  when $(u,x) \in   \partial B((v,A_j), \vert  A_j-A'_j \vert  )$ and hence --by the maximum principle--, also  for $(u,x) \in   { \mathcal C}_Y( \Sigma )\setminus B((v,A_j), \vert  A_j-A'_j \vert  )$. Taking $(u,x)=(0,x_0)$ we see that 
$c\;{ \frac {K_{(v,y)}(v,A_j)  } {K_{(v,A_j)} (v,A'_j)}} \leq 1$. So it follows from (\ref{HJB}) that 
\begin{align}\label{hbjohncy} K_{(v,y)}(u,x) \leq C'\, \sum _{j=1}^N \;\;K_{(v,A_j)}(u,x) \end{align} 
when $(v,y) \in   { \mathbb  R}\times M$, $ \vert  v \vert   \geq 1$,  and $(u,x) \in   { \mathcal C}_Y( \Sigma )$ satisfies $ \vert  u-v \vert   \geq 1$. Since by Theorem \ref{theorem41} $K_{(v,A_j)} \to  F_+$ for $j \to   \infty $ and since $F_+$ is minimal, the result follows. $\square$

The proof also yields the  following more general statement. Here the results of \cite{anc4} for John subsets with more than one pole  are effectively used.

{\bf Theorem \ref{presdubord}$\bm {}'$} {\it Let $ \{ M_n \}$ be a sequence of closed John regular subsets of $ \Sigma $ with  a common John constant $c_0$. Assume that $v_n \to  + \infty $ in ${ \mathbb  R}$ and that $(v_n, A_j^{(n)}) \to   \zeta  _j \in   {\bm  \Delta  }_1$, $1 \leq j \leq N$, where for each $n$,  $ \{ A_j^{(n)} \}_{1 \leq j \leq N}$ is a $c_0$-admissible set of poles for $M_n$. Then if $y_n \in   M_n$,  every Martin cluster function  of $ \{ (v_n,y_n) \}$ is   a linear combination of the $K_ {  \zeta  _j }$.}

In particular if  the points  $ \zeta  _j $ all coincide with a minimal boundary point $ \zeta  $ then $(v_n,y_n) \to   \zeta  $.

\subsection{\normalsize John cuts}\label{sectjohncuts}    Assume now that $M$ is a John regular closed subset of $ \Sigma $ and that $ \Sigma \setminus M$ is the disjoint union of two open subsets  $U_0$ and $U_1$. Fix $ \delta >0$ and  set $U_j^ \delta = \{ x \in   U_j\,;\, d_a(x,M) \geq  \delta \,  \}$ and $V_j^ \delta ={ \mathbb R} \times U_j^ \delta  $, $j=1,\,2$. Let $V_j={ \mathbb R} \times U_j$.

\begin{proposition}\label{newteo2}  If $h=K_\mu  $ is the positive superharmonic function in ${ \mathcal C}_Y( \Sigma )$ generated by a probability measure $\mu  $ supported on the closure of $V_0$ in the Martin compactification $ \widehat { \mathcal C}_Y( \Sigma )$ and not charging $(0, x_0)$, we have \begin{equation} \label{eq43} \hspace{2truecm}  h(u, x) \leq C\, [F_+(u,x)+F_-(u,x)],\hspace{0.7truecm} (u,x) \in   V_1^ \delta 
\end{equation}  for some constant $C=C( \Sigma ,M, U_0,x_0, \delta )$. 
\end{proposition}



{\bf Proof.}  We may assume that $x_0 \in   U_1^ \delta $ (using Harnack and changing the reference point) and  it suffices to prove (\ref{eq43}) for each $h=K_ {(v,y)} $, $y \in     V_0 $ with a constant $C >0$ as in the statement. Reducing $h$ on $ V_1 $,  it  suffices to prove the result for $K_{(v,y)}$, $y \in   M$, $v \in   { \mathbb  R}$.

For such a pole $(v,y)$, with say $v>0$, it follows from  (\ref{hbjohncy})  and   Theorem \ref{theorem41} that for $x \in   U_1^ \delta $ such that  $d(x,A_j) \geq {\frac  { 1} {2}} d(A_j,\partial  \Sigma )$ for $j=1,\dots ,N$, (we use the same  notations as above) 
\begin{align}K_{(v,y)}(u,x) &\leq C'\, \sum _{j=1}^N \;\;K_{(v,A_j)}(u,x) \nonumber\\
&\leq  C''\,\sum _{j} \,K_{(v,A_j)}(v,A'_j)\, e^{ \alpha _ {\rm max}  (u-v)}\,   \varphi  _0 (x) \nonumber \\ & 
\leq C'''\, e^{ \alpha_ {\rm max}  u}   \varphi  _0 (x).
\end{align} 
In the second line we have used the maximum principle (as above in the proof of Theorem \ref{presdubord}) to compare the positive $L$-harmonic function $e^{ \alpha_ {\rm max}  u}   \varphi  _0 (x)$ with the Green function with pole at $(v,A_j)$. In the last line we have used the inequalities given by Theorem \ref{theorem41} which imply that $K_{(v,A_j)}(v,A'_j)\simeq K_{(v,A_j)}(v-1,A_j)\simeq F_+(v-1,A_j)$. Using the similar inequality for $v \leq 0$ we get the desired conclusion. $\square$

\begin{remark} {\rm For $h=K_{(v,y)}$, $y \in   M$, the proof shows that  $h(u, x) \leq C\, [F_+(u,x)+F_-(u,x)]$ if $x \in     U_1$, $ \vert  u-v \vert   \geq 1$}. If moreover $v \geq 0$, then  $h(u,x) \leq C\, F_{+}(u,x)$, for these points $(u,x)$.
\end{remark}

\begin{corollary} \label{separation} Under the assumptions of Proposition \ref{newteo2} two sequences $ \{ (v_j,y_j) \}$ and $ \{ (w_j,z_j) \}$ such that $v_j,\, w_j \to   +\infty $, $y_j \in    U_0$ and $z_j \in    U_1 $, have at most one common cluster point in $ \bm  \Delta  $ which   can only be $F_+$. 
\end{corollary}


\begin{corollary} \label{} If $h_0$ and $h_1$ are positive harmonic in ${ \mathcal C}_Y( \Sigma )$, if $h_j=K_{\mu  _j}$ with $\mu  _j$ supported by $\overline {U_j \times { \mathbb  R}}$ (closure in $ \widehat { \mathcal C}_Y( \Sigma )$), $j=0,\,1$,  if $ \displaystyle \lim_{u  \to  + \infty } e^ {- \alpha _{max} u} (h_0\wedge h_1)(u,x_0)=\lim_{u  \to  + \infty } e^ {\alpha _{min} u}(h_0\wedge h_1)(-u,x_0)=0 $ then $h_0 \wedge h_1$ is a potential (that is, has no positive $L$-harmonic minorant).
\end{corollary}

\subsection {\normalsize Inner ball property} \label{sectinnerball}

If we have a boundary point $z \in   \partial  \Sigma $ and an open ball (or cap)  $B(a,r) \subset   \Sigma $ with $z \in   \partial B(a,r) $, $r<2$, the results in  \cite{anc2} (see also  \cite{anc0}) tell us (using Proposition \ref{classif}) that {\sl as  $v \to  v_0$ in ${ \mathbb  R}$ and $y \to  z$ non-tangentially in $B(a,r)$, the point $(v,y)$ tends to a finite type minimal boundary point $ \xi = \xi  (v_0;(z,a)) $ in  $ \widehat { \mathcal C}_Y( \Sigma )$.} Moreover the  minimal $K_ \xi  $ is bounded away from  $(v_0,z)$ and vanishes on $\partial  { \mathcal C}_Y(\Sigma )\setminus  \{(v_0,z)\} $. 

There is a parallel statement for  the behavior of $(v,y)$ for $v \to  + \infty $. But here the inner ball should be large. This will be used later for an example's construction in section \ref{sectionhemisphere}.

\begin{theorem} \label{theoreminnerhemisphere} Assume that $ \Sigma $ contains an open hemisphere $ \Sigma _+$ in $S_{d-1}$. Then if $ \{ y_j \}$ is a sequence in $ \Sigma ^ + $ and if $v_j \to  + \infty $ in ${ \mathbb R} $, the sequence $(v_j,y_j)$  converges in the Martin compactification of ${ \mathcal C}_Y( \Sigma )$ to the canonical Martin point $F_+$, i.e., $\lim _jK_{(v_j,y_j)}=F_+$.
\end{theorem}

The proof is deferred to  section \ref{sectionhemisphere}.

\section{Uniqueness of the Martin point at infinity for $d=3$.  }\label{sectdim3}

We now prove Theorem \ref{casdleq3co} (rather its cylinder version), using in an essential way a result of R.\ Ba\~{n}uelos and B.\ Davis on the heat kernels in planar domains (\cite{bada}, \cite{banProc}). This result says that given the point $x_1 \in    \Sigma $  there is a $t_0>0$ and for each $t \geq t_0$ a constant $C(t)>1$ such that  $\lim_{t\uparrow  \infty } C_t=1$ and 
\begin{equation}\label{BanBurEq} C(t)^{-1}\, e^{- \lambda _1t} \, \varphi  _0(x_1)\,  \varphi  _0(y) \leq   \pi  _t(x_1,y) \leq C(t)\, e^{- \lambda _1t} \, \varphi  _0(x_1)\,  \varphi  _0(y)
\end{equation}
when $t \geq t_0$ and  $y \in    \Sigma $ (see in \cite{bada} Theorem 1 and section 4). Recall $ \{  \pi  _t \}$ is the heat semi-group generated by the Laplacian in $ \Sigma $ and  $ \varphi  _0$ is normalized by the condition $ \Vert  \varphi  _0 \Vert _{L^2( \Sigma )}=1$.

\begin{theorem}\label{theoremd=3} If $d=3$, every sequence $  \{ \xi  _j \}_{j \geq 1}= \{ (v_j,y_j) \}_{j \geq 1}$,  in ${ \mathcal C}_Y( \Sigma )$ such that $v_j \to  + \infty $ converges to $F_+$, i.e., $K_ {\xi  _j}(u,x) \to  K_ { \zeta  _\infty } (u,x):= e^{ \alpha _{ \rm max } \,u}  \;\varphi  _0(x)/ \varphi _0(x_0)$ for $j  \to     \infty $. 
\end{theorem}

The following simple lemma (valid for all $d \geq 2$) deals  with times in $(0,t_0]$.

\begin{lemma} \label{lemma4.2} Given $ \delta _0>0$ and $x_1 \in    \Sigma$,  there is a constant  $C=C( \delta _0; \Sigma  ,x_1) \geq 1$ such that \begin{equation} \label{eq4.2}
 \pi_t(x_1,y)  \leq C e^{- \lambda_1 t} \varphi  _0(x_1) \varphi  _0(y)
\end{equation}
for all $y \in    \Sigma $ such that $ \vert  y-x_1 \vert   \geq  \delta _0$ and all $t \geq 0$
\end{lemma}

Assuming as we may that $ \delta _0< d(x_0,S_{d-1}\setminus  \Sigma )$, this is a simple consequence of the parabolic maximum principle applied in the region $ \{ (x,t)\,;\, t>0, \, x \in    \Sigma \, \}\setminus  \{ (x,t)\,;\,  \vert  x-x_0 \vert   \leq  \delta_0 , \, 0 \leq t \leq 1 \}$ (the two members of (\ref{eq4.2}) are $(\partial _t-L)-$harmonic in $(y,t)$ and the first has by definition minimal growth at infinity in $\Sigma\times  { \mathbb  R}_+   $ ).

\begin{lemma} Assume $d=3$ and let $T$, $\mu  >0$, $x_1 \in    \Sigma $ be given. Then, as $a \to  + \infty $,
\begin{equation}\label{removed} \int _0^T e^{-{ \frac { a^2  } {4t}}-\mu   t}  \,\pi  _t(x_1,y)\, { \frac { dt  } {\sqrt t}} =o( \int _T^ \infty  e^{-{ \frac { a^2  } {4t}}-\mu  t}  \,\pi  _t(x_1,y)\, { \frac { dt  } {\sqrt t}} )
\end{equation}
uniformly in $y \in    \Sigma $, $ \Vert y-x_1 \Vert  \geq  \delta _0$.
\end{lemma}

{\sl Proof.} We may assume $T \geq t_0$. By  (\ref{BanBurEq})  and lemma \ref{lemma4.2}, it suffices to prove the relation obtained from (\ref{removed})   when the terms $ \pi  _t(x_1,y)$ are removed from the integrals.

Now  $\int _0^T e^{-{ \frac { a^2  } {4t}}-\mu   t}  \;{ \frac { dt  } {\sqrt t}}  \;\leq  \;{ \frac { 1 } {2}}\,\sqrt T\;e^{-{ \frac { a^2  } {4T}}}$ and,  if $ \theta   \geq 1$,
$\int _ \theta  ^ \infty  e^{-{ \frac { a^2  } {4t}}-\mu  t}  \, { \frac { dt  } {\sqrt t}} \geq  e^{-{ \frac { a^2  } {4 \theta  }} }\; \int _ \theta  ^ \infty e^{-(\mu  +{ \frac { 1  } {2}})t}\, dt$, where for the last inequality we use the observation  that $t \mapsto { \frac { e^{ \frac { t  } {2}}  } {\sqrt t}}$ is increasing in $(1, \infty )$.

So $\int _ \theta  ^ \infty  e^{-{ \frac { a^2  } {4t}}-\mu  t}  \, { \frac { dt  } {\sqrt t}} \geq
{ \frac { 1  } {\nu }}\; e^{-{ \frac { a^2  } {4 \theta  }} }\;  e^{-\nu \theta  }$, where $\nu =\mu  +{ \frac { 1  } {2}}$, and we may conclude since as $a \to  + \infty  $, ${ \frac {  \nu  } {2}}\;\sqrt T\;e^{-{ \frac { a^2  } {4}}({ \frac { 1  } {T}}-{ \frac { 1  } { \theta  }})+\nu  \theta  } \to  0
$ for any fixed   $ \theta >T $. $\square$
\medskip

{\bf Proof of Theorem \ref{theoremd=3}.} Assume as we may that $y_j \to  y_ \infty  \in   \overline  \Sigma $.   Using the Ba\~nuelos-Davis Theorem and the above lemma and its proof, we have if $y_ \infty \ne x_1$ and $j \to   \infty $,
\begin{eqnarray}
\sqrt{4 \pi  }\,G(u,x_1; v_j,y_j) &=& \int _0^ \infty e^{-{ \frac { (u+(d-2)t-v_j)^2  } {4t}}}  \pi  _t(x_1,y_j)\; { \frac { dt  } {\sqrt t}}\cr
 &=&   e^{{ \frac { (d-2)(v_j-u)  } {2}}} \int _0^ \infty e^{-{ \frac { (u-v_j)^2  } {4t}}} e^{-{ \frac { (d-2)^2  } {4}}t} \;\pi  _t(x_1,y_j)\; { \frac { dt  } {\sqrt t}}\cr
\sim & e^{{ \frac { (d-2)(v_j-u)  } {2}}}&  \varphi _0 (x_1) \varphi_0  (y_j) \; \int _0^ \infty e^{-{ \frac { (u-v_j)^2  } {4t}}-[{ \frac { (d-2)^2  } {4}}+ \lambda _1]t } \; { \frac { dt  } {\sqrt t}}.
\end{eqnarray}
 Thus, for $x_1,\, x_2 \in    \Sigma $, we see that  $G(u,x_2; v_j,y_j)/G(u,x_1; v_j,y_j) \to   \varphi _0(x_2)/ \varphi  _0(x_1)$ for $j \to   \infty $ (assuming first   $y_ \infty \ne x_1$, $y_ \infty \ne x_2$). This  shows that a cluster function $K$ of  the Martin kernels $K_{v_j,y_j}$, as $j \to   \infty $,  is in the form $ K(u,x)=g(u)\varphi  _0(x)$ and hence must be $K=F_+$.

In fact, 
\begin{equation}\label{equivid}
\int _0^ \infty e^{-{ \frac { (u-v_j)^2  } {4t}}-[{ \frac { (d-2)^2  } {4}}+ \lambda _1]t } \; { \frac { dt  } {\sqrt {4 \pi  t\,}}}=   { \frac {  1   } {2 \sqrt{\mu \,}  }}\;  \,  e^{-\sqrt{A_j\mu  }}
\end{equation}
where $A_j= \vert  u-v_j \vert  ^2$, $\mu  ={ \frac { (d-2)^2  } {4}}+ \lambda _1$ (note that the left member of (\ref{equivid}) is the Green's function in ${ \mathbb R} $ with pole at the origin  and with respect to $ { \mathcal L}={ \frac {  d^2} {dx^2 }} -\mu  $ evaluated at $\sqrt {A_j}$). It follows that 
$G(u,x; v_j,y_j)/G(0,x_0; v_j,y_j) \to   F_+(u,x)$ when $j \to   \infty $.

This shows that $G((u,x_1);(v,y))\sim  {\frac  { 1} {\sqrt{(d-2)^2+4 \lambda _1}}} \varphi  _0(x_1)  \varphi  _0(y) e^ {- \alpha _{\max}(v-u)}$  uniformly with respect to $y$, as $v-u \to  + \infty $. Similarly, $G((u,x_1);(v,y))\sim  {\frac  { 1} {\sqrt{(d-2)^2+4 \lambda _1}}} \varphi  _0(x_1)  \varphi  _0(y) e^ {- \alpha _{\min}(u-v)}$    as $v-u \to  - \infty $,   uniformly w.\ r.\ to  $y$. $\square$ 
\vspace{6truemm}

{\bf Proof of Theorem \ref{theorem41}.} The proof of Theorem \ref{theorem41} (where $d \geq 2$) is completely similar. In fact as well known (see e.g.\  \cite{chakar}, \cite{pincho}) for every compact $K \subset   \Sigma $, $ e^{ \lambda _1t} \pi  _t(x,y) \to   \varphi  _0(x) \varphi  _0(y)$ as $t \to  + \infty $, uniformly w.\ r.\ to $y \in     K$. So an obvious adaptation of the above gives the convergence to $F_+$. Inequalities \ref{eq41} follow from the fact that now  $G(u,x_0;0,x_0) \sim C'\,e^{ \alpha _{\rm min} u}$ for $u \to   +\infty $.

\begin{remark} {\rm The proof also shows that Theorem \ref{theoremd=3} extends to $d \geq 3$ if the base $ \Sigma $ is intrisically ultracontractive with respect to $ \Delta _ S$, (see \cite{DS}, \cite{banProc}, \cite{ban91}). For example, by Ba\~{n}uelos results in \cite{ban91} this is the case if for some $p>{ \frac {  n} {2 }}$ the base $ \Sigma $ is $L^p$-averaging (that is  $ \rho _ \Sigma  \in   L^p( \Sigma )$ where $ \rho _ \Sigma (x)$ is the pseudo-hyperbolic distance to $x_0$). See \cite{ban91} for other examples.}
\end{remark}

\medskip

\section{Examples for $d  \geq 4$}\label{secexample1}

In this section we show that for $d \geq 4$ there are cones  in ${ \mathbb  R}^d$ with a host of   Martin points at infinity. 
See  sections \ref{example1construction} and \ref{example1description} (another example  described in section \ref{sectionhemisphere} shows that  these   points can be minimal as well as non minimal). This  is  closely connected with the existence --established by Cranston  and McConnell in \cite{CM}--  of a  bounded domain $D$ in ${ \mathbb  R}^3$ with an $h$-Brownian motion in $D$ with  an infinite  lifetime.  In fact we  use a  variant of the construction   in \cite{CM} section 3. 

  As before we  work with the model of the cone ${ \mathcal C}_0( \Sigma )$  given by the cylinder $({ \mathcal C}_Y( \Sigma ),L)$.

\vspace{4truemm} 

\subsection {\normalsize  Preliminary lemmas} \label{prell}

Fix $ \lambda_0> 0 $, $d \geq 2$, and consider a  cap $B=B(a,r):= \{ x \in   S\,;\,  \vert  x-a \vert  <r \}$  in $S:=S_{d-1}$, $r \leq 1$, with two given points $ \xi$, $ \xi  '   \in   \partial B$, symmetric in $S$ with respect to $a$. Let $T=B( \xi  ;r/100) \cap    \partial B$, $T'=B( \xi  ',r/100) \cap    \partial B$, and let $ M_B:= \{ x \in   B\,;\,  \Vert x- \xi   \Vert =  \Vert x- \xi  ' \Vert \,  \}$. 

Let $\Omega $ be a region in $S_{d-1}$ such that  $B \subset   \Omega  \subset   S_{d-1}\setminus (\partial B\setminus T)$, $\overline  {T'} \cap \overline {\Omega \setminus B}= \emptyset  $. Set $ \tilde \Omega =\Omega \times { \mathbb  R}$.

\begin{lemma}\label{lemma64}  Let $v=H_f $ solves  $ \Delta _S v- \lambda v-\partial _tv =0$ in $ \tilde \Omega  $ and $v(y,t)=\,f(y,t)$ in $\partial  \tilde \Omega $ where  $f(y,t)$ is bounded measurable   in $\partial  \tilde \Omega $, nondecreasing in $t$ and $f(y,t)=0$ for $y \notin   T'$.  Then, given $ \eta  \in   (0,1)$, 
  there exists  $   \varepsilon _1= \varepsilon _1(d, \lambda _0,  \eta )   >0$  such that for $0< \varepsilon  \leq  \varepsilon _1$ and  $0 \leq  \lambda  \leq  \lambda _0$,
\begin{eqnarray}\label{eqHf} H_f(x,t) \leq    \int_{\partial \Omega }\,     \, [\eta \,f(y,t)+(1-   \eta )f(y,t- \varepsilon r^2)] \;d\mu^\Omega   _x(y),\; (x,t) \in   M_B\times { \mathbb  R}\end{eqnarray}
Moreover $H_f(x,t)$ is nondecreasing in $t$. Here $\mu  _x^\Omega $ is the harmonic measure of $x$ in $\Omega $ w.\  r.\ to $ \Delta  _{S}- \lambda I$.
\end{lemma}

{\bf Proof.} The last claim follows from the parabolic maximum principle and the translation invariance in $t$ of  $ \Delta _S - \lambda I-\partial _t $.

To prove the first, observe that by the monotonicity assumption we may assume that  $f(y,s)= \varphi  (y)$ for $s>t- \varepsilon r^2$   and $f(y,s)=  \psi   (y)$ when $s \leq t- \varepsilon r^2$. Since the inequality is an identity when $f(y,s)$ is independent of $s$ we may assume $ \psi =0$ and also that $t= \varepsilon r^2$ by time translation invariance.

Then  $\int_{\partial \Omega }\,     \, [\eta \,f(y, \varepsilon r^2)+(1-   \eta )f(y,0)] \;d\mu^\Omega   _x(y)=       \eta  \Phi   (x)$ where $ \Phi $ solves $ \Delta   \Phi - \lambda  \Phi =0$ in $\Omega $  and  $\Phi =f$ on $\partial \Omega $. We want to show that 
$H_f(x, \varepsilon r^2)\leq   \eta    \Phi    (x)$ provided $ \varepsilon<{ \frac {  1} {8 }} $ is sufficiently small. Let $N$ be the integer part of ${ \frac {  4} { \varepsilon  }}$ and set $w(x,s)=H_f(x,s+ \varepsilon r^2)-H_f(x,s)$.

 By the  parabolic Harnack inequalities \cite{moser}, $Cw(a, r^2(1-k \varepsilon) ) \geq \,w(a, { \frac {  r^2} {2 }})$ for $1 \leq k \leq N$ and a constant $C=C(d)$. Thus, on summing up, $Nw(a,  r^2/2) \leq CH_f(a, r^2) \leq C \Phi (a)$. 

Applying next  the parabolic boundary Harnack principle in  $\Omega \times { \mathbb  R} $ (\cite{FGS}, \cite{heu}) to $w(x,s)$ and $ \Phi $ (viewed as functions of $(x,s) $) we obtain  $w( x,  \varepsilon r^2) \leq c_1\, {\frac  {w(a,r^2/2) } { \Phi (a)}}  \Phi (x) \leq { \frac {  c_1\, C} {N }}\,  \Phi (x)$ for $x \in   M_B$ with $c_1=c_1(C,d, \lambda _0)$. The result follows.
 $\square$

We will use lemma \ref{lemma64} in conjunction with the following  lemma.

\begin{lemma} \label{lemma65}Let $ \{  m _k \}_{1 \leq k \leq N}$ be a finite sequence of probability measures in ${ \mathbb  R}$ of the form $ m _k= {{\frac 1 2}}  \delta _0+{\frac 1 2} \delta _{-a_k}$ for $1 \leq k \leq N$, where $0 \leq a_k  \leq 1$. Let $L$ and $ \varepsilon$ be given positive numbers.
There is an $A=A(L, \varepsilon )>0$ such that if $  \sum_{k=1}^N a_k \geq A$, the measure $\nu_N  = m _1\ast \dots  \ast m _N $ satisfies : $ \nu_N  ([-L,0]) \leq  \varepsilon$.
\end{lemma}

The probability $ \nu  _N$ is the law of the random variable $Z:=- \sum _{j=1}^N a_j X_j$ if $ X_1,\dots, X_N$ are independent random variables such that $P(X_j=0)=P(X_j=1)={\frac 1 2}$. For $ \beta >0$, we have 
\begin{eqnarray}
\hspace{-0mm}P(-Z \leq L)=P(e^{ \beta Z}  \geq e^{- \beta L})  \leq e^{ \beta L} E(e^{ \beta Z})=e^{ \beta L} \prod_{j=1}^N E(e^{- \beta a_jX_j})=e^{ \beta L} \prod_{j=1}^N (1-{ \frac { 1-e^{- \beta a_j} } {2}}).\nonumber 
\end{eqnarray}
Thus, using $a_k \leq 1$, 
 $\displaystyle P(Z \geq -L) \leq e^{ \beta L}\,\prod _{j=1} ^N(1- { \frac {   \beta } {2 }}\,e^{-\beta }\, a_j) \leq e^{ \beta L}\, \exp(-\frac {   \beta } {2 } e^{- \beta }  \sum _{j=1}^Na_j\,)$. The lemma follows.\nobreak $\,\square$

\subsection {\normalsize  A class of cylinders.}\label{example1construction}

We now consider domains $ \Sigma  \subset   S_{d-1}$, $d \geq 4$, similar to  examples introduced  in  \cite{CM}: there are   disjoint open balls $B_j=B(x_j,r_j)$, $j \geq 0$, in $S_{d-1}$ such that (i) $ \sum_{j \geq 0} r_j^2=+ \infty $,  (ii) $B_j \subset    \Sigma $, (iii) for  $N \geq 1$, $ \Sigma  \setminus B_N$ has two components $ \Sigma  _N^+$, $ \Sigma  _N^-$ with disjoint closures and   $  \Sigma _N ^-  \supset   \bigcup_{j < N}    B_j$,   $  \Sigma _N^+\supset \bigcup     _{j>N} B_j$, (iv) there are   caps $ T_j$, $T'_j$ in $\partial B_j$, $j \geq 0$, symmetric with respect to $x_j$, of  radius $ \rho _j \leq r_j/ 10$ and such that $ \Sigma _j^- \cap    \overline  B_j \subset   T_j$, $ \Sigma _j^+ \cap    \overline  B_j \subset   T'_j$.

\begin{remark} \label{rksigma}{\rm There is an $ \varepsilon _d>0$ such that whenever $ r_j>0$, $j \geq 1$, satisfy   $ \sum r_j^2= \infty $ and $   \sum r_j^{d-1}   \leq  \varepsilon _d\,$, there exists a corresponding   $ \Sigma  $ such that moreover: (a)  $ \vert  x_j-x_k \vert   \geq 4 \, \max \{ r_k, r_j \}$ for $j\ne k$,  (b) the centers $x_j$ have a limit  $P_0$   in $S_{d-1}$ (c) $ \Sigma $ is  locally Lipschitz in $S_{d-1}\setminus  \{ P_0\}$ and is Dirichlet-regular in $S_{d-1}$. The proof  is  left to the reader.}\end{remark}

Set  $  \Sigma   _N= \Sigma _N^- \cup B_{N} $ for $N \geq 1$ and fix      $ \lambda _0 > 0$.  Let $k$ be a bounded positive solution of $ \Delta  k- \lambda k=0 $ , $k=0$ in $\partial  \Sigma _N\setminus T'_{N}$, $0 \leq  \lambda  \leq  \lambda _0$.  For $ \ell  >0$, let $h=h_ \ell$ solves: $ \partial _th(x,t)-\Delta _x h(x,t) + \lambda h(x,t)=0$ in $ \Sigma _N \times { \mathbb R} $, $h(t,x)=\bm 1_{t \geq -\ell}\, k(y)$ in $ \partial  \Sigma _N\times { \mathbb R} $.

\begin{proposition} Let $\ell $ and $ \varepsilon $ be positive reals and let $x \in     \Sigma _q^- $, $1 \leq q<N$.  There is an integer $ N_ \Sigma ( q,\varepsilon , \ell, \lambda _0)$ such that whenever $N \geq N_ \Sigma ( q,\varepsilon , \ell, \lambda _0)$, 
\begin{equation} \label{eq63} h_\ell (x,0)  \leq  \varepsilon \, k(x).
\end{equation}  
    
\end{proposition}

{\bf Proof.} Fix  $  \eta   ={ \frac {  1} {2 }}$ and a corresponding $ \varepsilon _0>0$ as given by lemma \ref{lemma64}. Let $p  \in    \{ q,\, q+1,\dots \}$, let $f(y,t)$ be a bounded Borel function in $T'_{p}\times { \mathbb R}  $ which is non decreasing in $t$ and let $v=H_f$ denote the solution of $ \Delta  _Sv- \lambda v-\partial _tv=0$ in $ \Sigma _p\times { \mathbb  R}$ with  $v=\bm 1_{T'_p\times { \mathbb  R}} \,f$ on the boundary. We show by induction on $n =p-q$,   that 
\begin{equation}\label{eq54} H_f(x,t)  \leq \int_{T'_p}  \, (\,\int _{- \infty }^ \infty \, f(y, t+s) \,d\nu  _{q,p}\,(s)\,)\, d\mu ^{ \Sigma _p} _x(y), \;\; x \in    \Sigma _q\setminus B_q 
\end{equation} 
where $   \nu  _{q,p}={\overset {p} {\underset {j=q} {\star } }} (  {\frac  { 1} {2}}   \, \delta _0+{\frac  { 1} {2}} \, \delta _{- \varepsilon _0r_j^2})$ and where $\mu^{ \Sigma _p}  _x$ is the harmonic measure of $x$ in $ \Sigma _p$ w.r.\ to $ \Delta  _ {S} - \lambda I$.  Denote $\nu _j:= {\frac  { 1} {2}}   \, \delta _0+{\frac  { 1} {2}} \, \delta _{- \varepsilon _0r_j^2}\,$.

For $n=0$ this is lemma \ref{lemma64}. Assuming that the property holds for $n-1 \geq 0$ and viewing  $H_f$ as a solution in $ \Sigma _{q}\times { \mathbb  R}$ of a Dirichlet problem for $  \Delta    _{S}- \lambda I-\partial _t$ we get by lemma \ref{lemma64} and maximum principle
\begin{eqnarray}
H_f(x,t) &\leq &   \int_{T'_q}   (\int_{- \infty }^{+ \infty }  \,H_f(y,s+t)\,  d\nu _q(s) \,)\;d\mu^{ \Sigma _q}   _x(y) \nonumber
\\ & \leq & \int_{T'_q}   (\int _{- \infty }^ {+\infty }  \,[\int_{T'_p}  \int _{- \infty }^ {+\infty }f(z, s+t+\tau )\, d\nu  _{q+1,p}( \tau )\, d\mu  _y^{ \Sigma _p}(z)\,  ]\, d\nu _q(s) \,)\;d\mu^{ \Sigma _q}   _x(y) \nonumber
\\ & = & \int_{T'_q}   (\int_{T'_p}  \,[ \int _{- \infty }^ {+\infty }  \int _{- \infty }^ {+\infty } f(z, s+t+\tau )\, d\nu  _{q+1,p}( \tau )\, d\nu _q(s) \,  ]\, d\mu  _y^{ \Sigma _p}(z) \,)\;d\mu^{ \Sigma _q}   _x(y) 
\nonumber 
\\ &=& \int_{T'_q}   (\int_{T'_p}  \,[ \int _{- \infty }^ {+\infty }   f(z, t+ \theta   )\, d\nu  _{q,p}(  \theta  )\,  ]\, d\mu  _y^{ \Sigma _p}(z) \,)\;d\mu^{ \Sigma _q}   _x(y) 
\nonumber
\\ &=& \int_{T'_p}   (\int _{- \infty }^ {+\infty }  f(z, t+ \theta   )\, d\nu  _{q,p}(  \theta  )\,  ]\, d\mu  _x^{ \Sigma _p}(z),\;\; x \in    \Sigma _q^-,  
\nonumber
\end{eqnarray}  
where we have used in the second line the induction assumption, in the third the fact that  integration with respect to $s$ and integration with respect to $z$ commute and --in the last line-- the  formula $\mu  _x^{ \Sigma _p}= \int \mu  _y^{ \Sigma _p}\, d\mu^{ \Sigma _q}  _x(y)$ (for $x \in    \Sigma _q^-$). This proves (\ref{eq54}).

From (\ref{eq54}) it follows that  for $x \in     \Sigma _q^-=\Sigma _q\setminus B_q$,
\begin{eqnarray} h_\ell(x,0) \leq  \int_{T'_N}  \, (\,\int _{- \ell }^ 0\, k(y) \,d\nu  _{q,N}\,(s)\,)\, d\mu ^{ \Sigma _N} _x(y) =k(x)\, \nu  _{q,N}([-\ell,0)), 
\end{eqnarray}
and the proposition follows from lemma \ref{lemma65} and the condition $ \sum_{j \geq 1} r_j^2 =+ \infty $. $\square$ 

We now take for $k$ the Green's function $k=G_y^ \lambda $ in $ \Sigma $ with pole at some point $y \in    \Sigma  \setminus \overline   \Sigma _N$ and with respect to $ \Delta_S  - \lambda I$. It is easily checked that for $-\ell  \leq s \leq 0$, $x \in    \Sigma _N$, $h_{\ell} (s,x) \geq  \int _{ 0 }^{\ell- \vert  s \vert  } e^{- \lambda t  }\, \pi _{t  }(x,y)\, dt$. Recall that $k(x)=\int_0^ \infty   e^{- \lambda t}  \pi  _t(x,y)\, dt$ and that the parabolic Green's function with pole at $(y_0,t_0 )$ in $ \Sigma $ --and w.r.\ to  $\Delta_{S}  - \lambda -\partial _t$-- is  $ \Gamma  (x,t;y_0, t_0):(x,t)\mapsto   \bm 1_{t>t_0} \, e^{- \lambda (t-t_0)} \pi  _{t-t_0}(x,y_0)$; thus $ u(x,s):=\int_{- \infty }^{+ \infty } \bm 1_{ \{ t>-\ell \}}\,  \Gamma  (x,s;y,t)\, dt= \int _{ 0 }^{\ell- \vert  s \vert  } e^{- \lambda t  }\, \pi _{t  }(x,y)\, dt$ is  bounded by $h_{\ell}(x,s)$ in $ \Sigma _N  \cap  \{ -\ell <s<0\}$ by the parabolic maximum principle. 

Thus, the previous result can be read  as follows.

\begin{lemma}\label{lemma67}    For any sequence $ \{ y_j \}$ converging in $ \Sigma  $ to the end ${ \mathcal E}$ of $\, \Sigma  $ defined by the cuts $B_N$ (a basis for ${ \mathcal E}$ is provided by the sets $\Sigma _N^+ $) and every $ \lambda  \geq 0$, we have
 $$\lim_{j \to   \infty } { \frac { \int _{0}^ {t_0}   e^{- \lambda t} \;\pi _s(x,y_j )\, ds } {\int _{0}^ \infty   e^{- \lambda t}\; \pi  _s(x,y_j )\, ds }}=0$$ for every fixed $t_0 \geq 0$ and every $x \in    \Sigma $.
\end{lemma}

So, $\lim_{j \to   \infty } { \frac { \int _{t_0}^ \infty   e^{- \lambda t} \pi _s(x,y_j )\, ds } {\int _{0}^ \infty   e^{- \lambda t} \pi  _s(x,y_j )\, ds }}=1$ for  $t_0 \geq 0$, $x \in    \Sigma $  and $ \lambda  \geq 0$.

\begin{lemma} \label{lemma68} Let $ \{ y_j \}$ be as in lemma \ref{lemma67}. For $j \to   \infty $ and for given   reals $ \rho $,  $\rho '$ the ratio
 \begin{equation}  {\frac { [\int_0^ \infty  t^{- {\frac { 1} { 2} }} e^{-{\frac {( \rho +(d-2)t)^2}{4t}}} \pi  _t(x,y_j)\, dt]} {[\int_0^ \infty  t^{- {\frac { 1} { 2} }} e^{-{\frac {( \rho' +(d-2)t)^2}{4t}}} \pi  _t(x,y_j)\, dt] } }
\end{equation}
converges towards $e^{- {\frac {d-2 } {2 } }( \rho - \rho ')}$ for each $x \in    \Sigma $.
\end{lemma}

We have $\int_0^ \infty  t^{- {\frac { 1} { 2} }} e^{-{\frac {( \rho +(d-2)t)^2}{4t}}} \pi  _t(x,y_j)\, dt=e^{{-\frac  { (d-2) \rho } {2}} }\int_0^ \infty  t^{- {\frac { 1} { 2} }}   e^{-{\frac  {  \rho ^2} {4t}}} e^{-t{\frac  { (d-2)^2} {4}} }\pi  _t(x,y_j)\, dt $. If we set $ \varphi  (t)=t^{-{\frac  { 1} {2}}  } e^{-{\frac  {  \rho ^2} {4t}}}e^{-t{\frac  { (d-2)^2}{4}}}$, then  for every $t_0>0$ and for $j \to   \infty $
\begin{equation} \label{eq67}\int _0^{t_0}  \varphi  (t)  \pi  _t(x,y_j)\, dt =o(\int _{t_0} ^ \infty  \varphi  (t)  \pi  _t(x,y_j)\, dt )
\end{equation}   
In fact, with $A>{\frac  { (d-2)^2} {4}} $, we have for $t_1>0$ large enough

\begin{equation}\int _{t_1} ^ \infty  \varphi  (t)  \pi  _t(x,y_j)\, dt  \geq C(t_1, A , d ) \int _{t_1} ^ \infty  e^{- A t}  \pi  _t(x,y_j)\, dt  
\end{equation}  (note that ${\frac  {\varphi  (t)  } { \varphi  (t_1)}}  \geq {\frac  {e^{- A t} } {e^{- A t_1}}}  $ for $t \geq t_1$, because $ \varphi  (t)e^{tA}$ is increasing for $t$ large enough). On the other hand for such a fixed $t_1$, we have  $\int _0^{t_1} \varphi  (t)  \pi  _t(x,y_j)\, dt  \leq C'(t_1, A,  \rho , d ) $ $\int _0^{t_1} e^{-At}\, \pi  _t(x,y_j)\, dt $ and (\ref{eq67}) follows by lemma \ref{lemma67}.

Since $e^{-{\frac  {  \rho ^2} {4t}} }  \to  1$ for $t \to   \infty $ ($ \rho $ being fixed) we see now that as $j \to   +\infty $,

\begin{equation} \int_0^ \infty  t^{- {\frac { 1} { 2} }}   e^{-{\frac  {  \rho ^2} {4t}}} e^{-t{\frac  { (d-2)^2} {4}} }\pi  _t(x,y_j)\, dt \sim \int_0^ \infty  t^{- {\frac { 1} { 2} }}    e^{-t{\frac  { (d-2)^2} {4}} }\pi  _t(x,y_j)\, dt
\end{equation}
Using also this  result for $ \rho ' $ the lemma follows. $\square$

\

\subsection{\normalsize The Martin boundary of the first example}\label{example1description}

 Using Lemma \ref{lemma68}  we  get  a (partial) description of the Martin boundary of the cylinder ${ \mathcal C}_Y( \Sigma )= { \mathbb  R}\times \Sigma  $ with respect to $ L:=\partial ^2_{uu} +{ \frac { d-2  } {2}} \partial _u+\Delta  _{S}$. In particular it will be seen that there are Martin boundary points related to sequences $(u_j,y_j) \in   { \mathcal C}_Y( \Sigma )$ with $ \lim u_j=+ \infty $ and distinct from  the canonical point $F_+$ given by  Theorem \ref{theorem2}.  

Denote $ \widetilde \partial  \Sigma :=  \partial  \Sigma  \cap     \{ \bigcup    _{n \geq 1} \partial  \Sigma _N^{-}  \}$ the set of points in $\partial  \Sigma $ ``away from the end ${ \mathcal E}$" (see Lemma \ref{lemma67}). Because ${ \mathcal E}$ is defined by a ``smooth" system of cuts (the balls $B_N$ or  the mediators $M_N$  of $T_N$ and $T'_N$ in $B_N$) it follows from standard forms of the boundary Harnack principle (see e.g.\ \cite{anc2}, \cite{anc4})  that for $ \lambda < \lambda _1( \Sigma )$, the end ${ \mathcal E}$ is the trace on $ \Sigma $ of the neighborhoods system --in the Martin compactification of $( \Sigma , \Delta_S  + \lambda I)$-- of a minimal $ \Delta _S + \lambda I$--harmonic function $k_{ \mathcal E}^ \lambda $ (normalized at $x_0$) which   vanishes on $ \widetilde \partial  \Sigma $ (\cite{anc2} Th\'{e}or{\`e}me 2.5). 

For  $ \alpha  \in   [ \alpha _0, \alpha _ {\max })$ --recall $ \alpha _{0} :=-{ \frac { d-2  } {2}}  $  and $\alpha _{ \max } : = {\frac {-(d-2)+\sqrt{(d-2)^2+4 \lambda _1( \Sigma )} } {2 } }$--, define $K_{ \mathcal E}^\alpha (u,x):=e^{ \alpha u}\, k_{ \mathcal E}^ {\lambda ( \alpha )}(x)$, $(u,x) \in    \Sigma \times { \mathbb  R}$. Here $ \lambda ( \alpha )= \alpha ^2+(d-2) \alpha $ (thus $ \lambda ( \alpha _0)=-{ \frac { (d-2)^2  } {4}}$). Recall $ {\bm \Delta }  $   denotes the Martin boundary of ${ \mathcal C}_Y( \Sigma )$ w.r.\ to $L$, and $ {\bm \Delta }_1$ its minimal part.

\begin{theorem}\label{theoremexample1} The function 	$K_{ \mathcal E}^{ \alpha _0}$ is $L$-minimal (so $K_{ \mathcal E}^{ \alpha _0}=K_ \xi  $ for some $ \xi   \in   {\bm  \Delta  }_1$) and there exists $ \Phi : \Sigma  \to  { \mathbb  R}_+$ going to $  + \infty $ along ${ \mathcal E}$ and such that $ (u_j,y_j)  \to   \xi $ when $ \vert  u_j \vert   \leq  \Phi (y_j)$ and  $ \{ y_j \} \to  { \mathcal E}$.  Moreover for every $ \alpha  \in   ( \alpha _0, \alpha _{ \max })$, $K_ { \mathcal E}^ \alpha $ is minimal $L$-harmonic in ${ \mathcal  C}_y( \Sigma )$ and each  associated sequence $ \{ u_j,y_j\}$ in ${ \mathcal C}_Y( \Sigma )$ satisfies : (i) $ \{ y_j \} \to  { \mathcal E}$ and (ii) $u_j \to  + \infty $. 
\end{theorem}

Similarly, for $ \alpha  \in   ( \alpha _{\min}, \alpha _0)$, the function $K_ { \mathcal E}^ \alpha $ is $L$-minimal in ${ \mathcal  C}_y( \Sigma )$ and every associated sequence $ \{ u_j,y_j\}$  satisfies : (i) $ \{ y_j \}$ converges to ${ \mathcal E}$ and (ii) $u_j \to  - \infty $.

{\bf Proof.} (a) By lemma \ref{lemma68}, if $ \{ y_j \}$ is a  sequence in $ \Sigma $ converging towards ${ \mathcal E}$ and   such that $ \{ (0,y_j) \}$ converges to some $    \xi   \in   {\bm \Delta }  $, the Martin function $K_ \xi  $ satisfies:  $K_ \xi  ( \rho ,x)/K_ \xi  ( \rho ',x)=e^{-{\frac  { d-2} {2}  (\rho - \rho ')}} $ for $ \rho ,\,  \rho ' \in   { \mathbb R} $ . Thus $K_ \xi  (u,x)=e^{- {\frac {d-2 } {2 } }u}\,k(x)$ where $k$ is independent of $ u $ and necessarily a positive solution of $ \Delta  _{S} k+ \lambda_0 k=0$ in $ \Sigma $, $ \lambda_0 = \lambda ( \alpha _0)$.

Using the John cuts ${M_n}$ (the mediator in $B_n$ between $T_n$ and $T'_n$), $n \geq 1$,  and Proposition \ref{newteo2} we see that $h_ \xi  $ vanishes on $ \tilde \partial  \Sigma \times { \mathbb R} $. So $k_ \xi  =0$ in $ \tilde \partial  \Sigma $ and as mentioned before $k$ must be the $( \Delta  _{S_{d-1}}+ \lambda _0I)$-minimal function corresponding to ${ \mathcal E}$, i.e $k= k_{ \mathcal E}^ {\lambda ( \alpha _0)}$.

(b) It  follows that $(u_j,y_j) \to    \xi  $ when $ \{ u_j \}$ is bounded and  $ \{ y_j \}$ as before (see Proposition \ref{classif}). And  for $ \rho _y$  growing sufficiently slowly to $+ \infty $  as $y \to  { \mathcal E}$,
the point $  ( u,y) $ tends to $\xi  $ for $y \to  { \mathcal E}$ and $ \vert  u \vert   \leq  \rho _y$ (the convergence holds in the Martin space of $({ \mathcal C}_Y( \Sigma ),  L)$\,). 
In particular there is no minimal boundary point $ \zeta  =\lim (u_j, y_j)$, with $y_j \to  { \mathcal E} $ of the finite type (i.e.\ non exponential in the vertical variable) described in Proposition \ref{classif} (iii).


(c) We now show that  $h_ \xi  $ is minimal $L$-harmonic in ${ \mathcal C}_Y( \Sigma )$ and more generally that for each $ \alpha  \in   [ \alpha _0, \alpha _ {\max} )$ the function   $h_0(u,x) =e^{ \alpha u}\, k_0(x)$, where $k_0=k_{ \mathcal E}^{ \lambda ( \alpha )}$,  is minimal harmonic for $({ \mathcal C}_Y( \Sigma ), L)$. Consider its Martin's disintegration into $L$-minimal functions. This disintegration does not charge the set of minimal functions in the form $e^{ \beta u}k(x)$ with $ \beta \ne \alpha $ (this would contradict the behavior of $h_0$ as $u\uparrow + \infty $ or $u\downarrow - \infty $). Since  $h_0=0$ on $ { \mathbb R}\times  \widetilde \partial   \Sigma   $, it  is supported by the set of  minimal points  $ \zeta  $ such that $  \zeta    =\lim ( u_j, y_j)$, $ \vert  u_j \vert   \to   \infty $ and $y_j \to  { \mathcal E}$ (using Proposition \ref{newteo2} and (b) above) and hence $K_ \zeta  (u,x)=e^{ \beta u}k(x)$ with $k$ vanishing on $ \tilde \partial  \Sigma $.  Thus $k=k_{ \mathcal E}^{ \lambda (  \beta  )}$ and the disintegration is supported by a minimal point $ \zeta  $ such that $K_ \zeta  =h_0$. 

The remaining assertion clearly follows from Corollary 
\ref{separation} and the proof is complete. $\square$


\vspace {1truecm} 

\subsection{\normalsize Proof of Theorem \ref{theoreminnerhemisphere} and a second  example}\label{sectionhemisphere}

 In this section we construct a second example --based on the first--  for which there are for each $ \alpha  \in   ( \alpha _{\min},  \alpha _{\max})$  corresponding minimal and non minimal Martin points $ \zeta  $ with  $K_{ \zeta  }(u,x)=e^{ \alpha u}\, k(x)$.
We first establish --in the spirit of \cite{anc0}-- Theorem \ref{theoreminnerhemisphere}.

{\bf A. Preliminaries.} Denote $ \Sigma  _+$  the hemisphere $ \{ t \in   S_{d-1};\,t_1>0 \}$ of $S_{d-1}$, $x_0=(1,0,\dots,0)$ its center, and   $ \sigma  $  the reflexion $x=(t_1,\dots, t_d) \mapsto (-t_1,t_2, \dots, t_d)$.

 \begin{proposition} \label{ineqfondhemis} Assume that $ \Sigma \supset  \Sigma _+$ and denote $G$ the Green's function of ${ \mathbb  R}\times  \Sigma $. Given $ \rho >0$, there is  a $C=C(d, \rho )$ such that whenever $y\in    \Sigma  _+$, $x  \in    \Sigma $, and $u \leq v-2 \rho $
\begin{align} G_{(v,y)}(u,x):=G(u,x;v,y) \leq C\, G_{(v,y)}(u,x_0)
\end{align}

\end{proposition} 


{\sl Proof.} a) By a known  general estimate (see \cite{anc0},\, \cite{anc2}) of the Green's function of a domain containing a $C^2$-ball (here $ \{ (t,x) \in   { \mathcal C}_Y( \Sigma )\,;\,  \vert  x-x_0 \vert  ^2+ \vert  t-w \vert  ^2 \leq 2 \}$)  -together with Harnack inequalities and elementary observations- we have when $z \in    \Sigma _+$, $z' \in    \Sigma$,  $\vert  w'-w  \vert  \geq  \rho $ :
\begin{align}\label{ineqcarbou}G(w',z'; w,z) \leq  C_{d, \rho }\, \,G(w\pm  \rho ,x_0;w,z)
\end{align}

b) Denote $ \tilde  \sigma (t,x) = ( t,\sigma  (x))$. Applying the  maximum principle in $ { \mathbb R} \times \Sigma _+$ to the functions $G_{(v,y)}$ and $G_{(v,y)} \circ   \tilde \sigma  $ in ${ \mathcal C}_Y(  \sigma  (\Sigma _+))$ -extend the second by zero outside ${ \mathcal C}_Y(  \sigma  ( \Sigma))$- we have
\begin{align} \label{ineqreflex} G_{(v,y)}(w,z) \leq  G_{(v ,y)}(w , \sigma  (z)) \leq \sup \{G _{(v ,y)}(w , z')\,;\, z' \in    \Sigma _+ \} 
\end{align}
for $z \in    \Sigma _-= \Sigma \setminus  \Sigma _+$, $w \ne v$. In particular we need only consider  $x \in    {\Sigma} _+  $  (take $w=u$). 

Similarly, if $\mu  _{(u,x)}^t$, $t \in   { \mathbb R} $,  is the $L$-harmonic measure of $(u,x)$, $u \leq t$,  in the truncated cylinder $C_-^t:={ \mathcal C}_Y( \Sigma ) \cap     \{ (w,z)\, ;\, w<t\,  \}$ we have for $x \in    \Sigma _+$
\begin{align}\label{ineqharmomeas} {\ } \hspace {35truemm} \mu  _{(u,x)}^t \leq   \sigma  (\mu  _{(u,x)}^t) \hspace{10truemm} {\rm \ in \ \ } \{ t \}\times   \sigma  (\Sigma _-) \subset   \partial C_-^t\end{align} 
This is because  the adjoint Green's function ${\stackrel * G}{}^{\,t}_{(u,x)}$ with respect to $C_-^t$ is larger than ${\stackrel * G}{}^{\,t}_{(u,x)} \circ   \tilde \sigma  $ in $C_-^ t \cap   \{ (w,z)\,;\, z \in     \Sigma _+ \}$, and $\mu  _{(u,x)}^t(t,dz)$ is  $\partial _w{\stackrel * G}{}^{\,t}_{(u,x)}\, d \sigma  _S(z)$.

c) Now write the r\'{e}duite (w.\ r.\ to  the cylinder ${ \mathcal C}_Y( \Sigma )$) of $G_{(v,y)}$    over $ \Sigma_+ ^{u+ \rho }:= \{ u+ \rho \}\times  \Sigma_+ $, i.e.,  $p=R_{G_{(v,y)}}^{ \Sigma_+ ^  {u+ \rho }} $ (\cite{brelot}, \cite{CC}),  as a potential $G_\mu   $ of a positive measure $\mu  $ on $ \Sigma_+ ^{u+ \rho }$. Then 
\begin{align}\label{ineqcar2}G_{\mu }(u,x)&=\int _{ \Sigma_+ ^{u+ \rho }}G(u,x; \zeta  )\, d  \mu ( \zeta  )  \nonumber 
\\ &\leq C\, \int _{ \Sigma _+^{u+ \rho }}G(u,x_0; \zeta  )\, d  \mu  ( \zeta  ) \nonumber
\\ & = C\, G_{\mu }(u,x_0) \leq C\, G_{(v,y)}(u,x_0), \end{align}
 using  (\ref{ineqcarbou}). 

d) Finally $ q  =G_{(v,y)}-p$ is  majorized in $C^{u+ \rho }_-$ by the solution $h$ to the Dirichlet problem in $C^{u+ \rho }_-$, with the boundary condition $ h  =G_{(v,y)}$ in $ \Sigma _-\times  \{ u+ \rho \,  \}$ and $ h  =0$ on the rest of $\partial U$. By (\ref{ineqreflex}) and (\ref{ineqharmomeas}),  $q(u,x) \leq p(u,x)$ and by (\ref{ineqcar2}) $q(u,x) \leq C\, G_{(v,y)}(u,x_0)$. $\square$


\begin{corollary} \label{corinnerhemisphere} Suppose $ \Sigma _+ \subset    \Sigma $. Let $ \{ y_j \}$ be a sequence in $ \Sigma ^ + $ and let $v_j \to  + \infty $ in ${ \mathbb R} $. Then $ \{ (v_j,y_j) \}$ converges to the canonical Martin point at $+ \infty $, i.e., $\lim _jK_{(v_j,y_j)}=F_+$.
\end{corollary}

This is Theorem \ref{theoreminnerhemisphere}. Similarly of course, if $v_j \to  - \infty $, $y_j \in    \Sigma _+$, we have $\lim _jK_{(v_j,y_j)}=F_-$.

{\sl Proof.}  If $ \zeta  $ is a cluster Martin point for $ \{ (v_j,y_j) \}_{j \geq 1}$ it immediately follows from   Proposition \ref{ineqfondhemis}  that if $h$ denotes the $L$-harmonic measure of $ \{ 0 \}\times  \Sigma $ in ${ \mathcal C}_Y( \Sigma )$, 
 \begin{align} K_ \zeta  (u,x) \leq C\, K_{ \zeta  }(u_0,x_0) \, h(u-u_0,x)  \nonumber
\end{align}
for $u \leq u_0$, $x \in    \Sigma $. This shows  that $K_ \zeta  $ vanishes on  ${ \mathbb R} \times \partial  \Sigma $ and that $K_ \zeta  $  is  bounded for $u \leq u_0$. Thus $K_ \zeta  =F_+$ by Theorem \ref{theorem2}. $\square$

\vspace{10truemm}
{\bf B.} Assuming $d \geq 4$ we construct  $ \Sigma $ as follows. We start with the hemisphere $ \Sigma _+$, a point $P_0 \in   \partial  \Sigma _+$  and a sequence of points  $P_n \in   \partial  \Sigma _+$ such that  $ \vert P_n-P_ 0 \vert  =4^{-n}$, $n \geq 1$. For each $n \geq 1$, let $\Omega _n$ be a domain in $S_{d-1}$, $\overline  \Omega _n \subset   B_{S_{d-1}}(P_n, {\frac  { 1} {4^{n+1}}} )  \setminus   \overline   \Sigma _+ $  of the type considered in  remark \ref{rksigma},  (starting with disjoint balls $\overline  B_{n,j}$, $j \geq 1$, in $S_{d-1}$   converging to some point $Q_n$ in $B(P_n, {\frac  { 1} {4^{n+1}}} )  \setminus   \overline   \Sigma _+$ and such that the sum of the squares of their radii diverges). Let $U_n$ be a  region in $B(P_n, {\frac  { 1} {4^{n+1}}} )  \setminus   (\overline   \Sigma _+\cup \overline  \Omega _n)$ such that $\overline  U_n  \cap    \overline   \Sigma _+$ is a closed ball $ \overline  \Delta  _n $ in $\partial  \Sigma _+$ of center $P_n$ and $\overline  U_n \cap    \overline  \Omega _n$  a  cap $ \overline  \Delta  '_n \subset   \partial B_{n,1} \cap    \partial \Omega _n$ (where $ \Delta  _n$ and $ \Delta  '_n$ are the relative interiors).

The domain  $\Sigma $  is  the union of $ \Sigma _  +$, the joining regions $U_n$ and the disks $ \Delta  _n$ and $ \Delta  '_n$, $n \geq 1$. 

Let $ \alpha  \in   ( \alpha _{\rm min} , \alpha _{\rm max })$ and $n \geq 1$. As before, by  \cite{anc2} there is a unique positive $(\Delta  + \lambda ( \alpha )I)$-harmonic function $k_n$ in $ \Sigma $ vanishing in $\partial  \Sigma \setminus  \{ Q_n \}$ and such that $k_n(x_0)=1$ ($x_0$ is  the center of the hemisphere $ \Sigma $). Moreover $k_n$ is minimal $ \Delta  $-harmonic in $ \Sigma $ and  by Theorem \ref{theoremexample1} the function $K_{Q_n}^{( \alpha )}(u,x)=e^{ \alpha u}\, k_n(x)$ is  minimal $L$-harmonic  in ${ \mathcal C}_Y( \Sigma )$. Denote $h^ {( \alpha )}(u,x)=e^{ \alpha u}\, k(x)$ the similar $L$-harmonic function in ${ \mathcal C}_Y( \Sigma )$ with a pole at $P_0$.

\begin{proposition} The function $h^{( \alpha )}$ is a non minimal Martin function for $({ \mathcal C}_Y( \Sigma ),L)$.
\end{proposition}

{\sl Proof.} It follows from standard arguments  that the $k_n$ vanish uniformly on the boundary of $ \Sigma $ away from $P_0$ (as $n \to   \infty $), so that $k=\lim k_n$. Thus $h^{( \alpha )}=\lim K_{Q_n}^{( \alpha )}$ and $h^{( \alpha )}$ is an $L$-Martin function in ${ \mathcal C}_Y( \Sigma )$ associated to a point $ \zeta  \in   \bm  \Delta  $. If $ \alpha  \geq  \alpha _0$ (resp. $ \alpha  \leq  \alpha _0$),  there is a sequence $ \{ (v_j,y_j) \} $ with $y_j \in   \Omega _j$, $\lim v_j=+ \infty $ (resp. $\lim v_j=- \infty $) converging to $ \zeta  $.

By Corollary \ref{corinnerhemisphere}, the point $ \zeta  $ --as a Martin boundary point-- is not in the closure of ${ \mathbb R} \times  \Sigma _+$. Thus every sufficiently small neighborhood $V$ of $ \zeta  $ meet ${ \mathbb R} \times \Omega _n$ for all large $n$, but not ${ \mathbb R} \times \Sigma_+ $. And $V \cap    { \mathcal C}_Y( \Sigma)$ is not connected. But  (by a general property) each neighborhood of   a {\sl minimal} Martin boundary point contains another whose trace in ${ \mathcal C}_Y( \Sigma )$ is connected (see e.g. \cite{naim} p.\ 223). Hence $ \zeta  $ is not minimal. $\square$

\section{ Extensions to more general cylinders}\label{sectgene}

The argument in sections 2 and 3 can be extended to more general second order elliptic operators in cylinders. We describe here a simple generalization and state the corresponding results. 
Assume $\Sigma $ is a {\em relatively compact} region in a $C^1$ Riemannian manifold $M$ of dimension $d-1$ ($d \geq 2$) equipped with a second order  elliptic operator $L_M$ in the form 
\begin{equation}L_M( \varphi  )= {\rm div}(A\,\nabla  \varphi  )+ B.\nabla  \varphi   + \gamma \, \varphi  
\end{equation}
where $A$ is a measurable, bounded  and uniformly elliptic section of ${\rm End}(T(M)$, $B$ a bounded measurable vector field in $M$ and $ \gamma $ a nonpositive bounded measurable function in $M$ (ref.\ \cite{sta}). We also assume  that $M\setminus  \Sigma $ is non polar. Thus $L_M$ admits a Green function in $ \Sigma $.

We consider now  a differential  operator in the cylinder ${ \mathcal C}_Y( \Sigma )= { \mathbb  R}\times  \Sigma $ which is a direct sum $L=L_{ \mathbb R} \oplus L_M$ where $L_{ \mathbb R}= {\frac { d^2} { du^2} } + b  {\frac {d } {du } } $ is  translation invariant  in ${ \mathbb R} $ (i.e., $b $ is a real constant).  
Again we fix  some $ x_0 \in    \Sigma $ and take $(0,x_0)$ as the normalization point in ${ \mathcal C}_Y( \Sigma )$ for the Martin functions. 

We define  $ \lambda _1:= \lambda _1(L_M; \Sigma ) $ as the supremum of all real $t$ such that  $L_M+tI$ admits a Green's function in $ \Sigma $  (or such that the cone of nonnegative $L_M+tI$-superharmonic functions has a dimension $>1$). It is well known  that $ 0<\lambda _1< \infty $ and that  for $t:= \lambda _1$ all nonnegative $(L_M+tI)$--superharmonic functions in $ \Sigma $ are proportional to the unique (up to scalar multiplication) positive $L_M+tI$ positive solution $ \varphi  _0$ in $ \Sigma $. This solution $ \varphi  _0 $ is bounded, vanishes in the weak sense on the boundary $\partial  \Sigma $ and $ \varphi  _0 \in   H_0^1( \Sigma )$. As also well-known we have similar properties for the formal adjoint operator $L_M^\ast$, and moreover $ \lambda _1(L^\ast_M; \Sigma )=\lambda _1(L_M; \Sigma ) $. We denote $  \varphi  _0^\ast$ a  positive $(L^\ast+ \lambda _1I)$-superharmonic function in $ \Sigma $. Again $ \varphi  ^\ast_0$ is unique up to multiplication by a constant, vanishes on $\partial  \Sigma $ and is $(L_M^\ast+ \lambda _1I)$-harmonic in $ \Sigma $.

Repeating the argument used in section \ref{secttheorem1} we obtain  a  similar description of the $L$-minimal Martin function associated with the end $u \to  + \infty $ in ${ \mathcal C}_Y( \Sigma )$.

\begin{proposition} \label{minimalinfinitygenral}
If $K$ is a minimal $L$-Martin function in ${ \mathcal C}_Y( \Sigma )={ \mathbb  R}\times  \Sigma $ associated to a sequence $(u_j,x_j)$ with  $u_j \to  + \infty $, $x_j \in    \Sigma $, then $K$ is in the form   
\begin{equation}K(u,x)=e^{ \alpha u} K(0,x),\;\; (u,x) \in   { \mathbb  R}\times    \Sigma    \end{equation}
for some $ \alpha  \geq  -{\frac { b} { 2} } $  and $s(x)=K(0,x)$ is a proper function in $ \Sigma $: $ L_M  s+ \lambda s=0$, $ \lambda = \alpha ^2+b \alpha $. So $ \alpha ={ \frac {-b+\sqrt{b^2+ 4\lambda }   } {2}} $ and $ -{ \frac { b^2  } {4}} \leq \lambda  \leq  \lambda _1$; moreover $s$ is $( L  _M+ \lambda I)$--minimal in $ \Sigma $. 
\end{proposition}

As before there is  a natural bijection $K \mapsto   \tilde K$  between the set ${\bm  \Delta  }^{+ \infty }$ of the Martin function arising from some sequence $(v_j,y_j)$ in ${ \mathcal C}_Y( \Sigma ) $ with $\lim v_j=+ \infty $ and the analogue set ${\bm  \Delta  }^{- \infty }$ (related to the condition $\lim v_j=- \infty $) by letting $ \tilde K(u,x)=e^{-bu} K(-u,x)$. 

Theorem \ref{theorem2bis} can also be extended  to the present framework, but a slight modification is required in the proof. Set $F_+(u,x)=e^{ \alpha _{\rm max} u}  {\frac  { \varphi  _0(x)} {\varphi  _0(x_0)}} $ and $F_-(u,x)=e^{ \alpha _{\rm min} u}  {\frac  { \varphi  _0(x)} {\varphi  _0(x_0)}}$ where $\alpha _{ \max } : = {\frac {-b+\sqrt{b^2+4 \lambda _1} } {2 } }$ and $\alpha _{ \min } : = {\frac {-b-\sqrt{b^2+4 \lambda _1} } {2 } }$.

\begin{theorem}\label{theorem2ter} If $f$ is $L$-harmonic in ${ \mathcal C}_Y( \Sigma )$ and vanishes in weak sense on ${ \mathbb  R}\times \partial  \Sigma $, then $f$ is a linear combination of $F_+$ and $F_-$. In particular, if moreover $ \displaystyle \liminf_{u \to  - \infty } e^{ \alpha _{\rm min} u} f(u,x_0)=0$ then $f$ is proportional to $F_+$. Thus $F_+$ and $F_-$ are $L$-minimal in ${ \mathcal C}_Y( \Sigma )$.
\end{theorem}

As in the proof of Theorem \ref{theorem2bis},
we may reduce ourselves  to show the following. A function $F$ in ${ \mathcal C}_Y( \Sigma )$ which  vanishes on ${ \mathbb R} \times \partial  \Sigma $ and which is in the form $F(u,x)=\int _{A'}\, K_ \zeta  (u,x)\, d \mu  (\zeta ) $ where $A'=  \{  \zeta   \in  {\bm  \Delta  }^{+ \infty } \cap    {\bm  \Delta  }_1\,;\,  -{\frac  { b} {2}} \leq  \alpha( \zeta  )  \leq  \alpha _1\} $, $ \alpha _1< \alpha _{\max} $,  and where $\mu  $ is  a finite positive Borel measure on $A$, must be the zero function. Denote $ \lambda '_1= \alpha_1 ^2+b \alpha _1$.

As before the function $ \varphi  (x)= \int_{A'} K_ \zeta  (0,x)\, d\mu  ( \zeta  )$ is positive superharmonic with respect to $L_0=L_M+\lambda _0I$, $ \lambda _0=-{\frac  { b^2} {4}} $ and vanishes in the weak sense on $\partial  \Sigma $. The measure $ -L_0( \varphi  )$ is given by the density $\psi (x)  =\int_{A}  (\lambda- \lambda _0)\, k_  \zeta    ^  \lambda  (x)\, d\mu  _f(x)$. It follows that $  \varphi  $ is the $ L_0 $-Green's potential in $ \Sigma $ of $ \psi $ and again, 
$ G_ \Sigma ^{ L_0 }( \psi )= \varphi   \geq   {\frac { 1} { \lambda '_1 - \lambda _0} }  \psi $
in $ \Sigma $, where $G^{ L_0  }_ \Sigma$ is  Green's function in $ \Sigma $ w.r.\  to $ L_0$. To  conclude we then slightly modify the argument in section \ref{secttheorem1} using now the minimal heat semi-group $P_t$ generated by $L$ in $ \Sigma $.
 \begin{align}   
 \int_ \Sigma  G_ \Sigma ^{L_0}( \psi )\,  \varphi   _0^\ast\, d \sigma  _M=\int _0^ {+ \infty }\int _ \Sigma \,e^{ \lambda _0t}   P_t( \psi )\,\varphi  _0^\ast\, d \sigma  _M\, dt=\int _0^ {+ \infty }\int _ \Sigma e^{-( \lambda _1- \lambda _0)t}  \psi \, \varphi  _0^*\, d \sigma  \, dt \nonumber
\end{align}
because $P_t^{\ast}( \varphi  _0)=e^{- \lambda _1t}  \varphi  _0^*$. Thus $\int_ \Sigma  G_ \Sigma ^{L_0}( \psi )\,  \varphi   _0^\ast\, d \sigma  _M={\frac  { 1} { \lambda _1- \lambda _0}} \, \int  \psi \,  \varphi  _0^*\, d \sigma  _M$. 

But on the other hand from  $ G_ \Sigma ^{ L_0 }( \psi ) \geq   {\frac { 1} { \lambda '_1 - \lambda _0} }  \psi $
it follows  that $\int_ \Sigma  G_ \Sigma ^{L_0}( \psi )\,  \varphi   _0^\ast\, d \sigma  _M$ is larger than ${\frac { 1} { \lambda '_1 - \lambda _0} } \int  \psi\,  \varphi  _0^*\, d \sigma   $. Thus $\int   \psi   \,  \varphi  _0^\ast \, d \sigma  =0$ and $  \psi  = 0$ in $ \Sigma $. So $ \varphi    =0$. $\square$ 
 
\begin{corollary}\label{theorem4} Every positive $L$-harmonic function $f(u,x)$ on ${ \mathcal C}_Y( \Sigma )$ vanishing (in the weak sense) on ${ \mathbb R} \times \partial  \Sigma $  and such that  $ \displaystyle \limsup_{u \to  - \infty}  f(u,x_0)< \infty  $ --$x_0 \in    \Sigma $-- is proportional to $F_+$.

\end{corollary}

 Theorems \ref{theorem41} and \ref{theoremd=3}   extend as follows.

\begin{theorem}  Let $ \xi  _j=(v_j,y_j)$, $j \geq 1$, be a  sequence of points   in ${ \mathcal C}_Y( \Sigma )$ such that $v_j \to   +\infty $. If $d \leq 3$ or if $d \geq 4$ and $ \{ y_j \}$ is relatively compact in $ \Sigma $, the functions $K_ {\xi  _j}(u,x)$ converge to $ K_ { \zeta  _\infty } (u,x):= e^{ \alpha _{\rm  \Sigma } \,u}  \; \psi   _0(x)$. In particular when $d \leq 3$, $ \zeta  _ \infty  $  is the only Martin point at infinity. 
\end{theorem}

The proofs are the same as above  in section \ref{sectdim3} using the natural extensions of (\ref{BanBurEq}) (with $y$ in a relatively compact subset of $ \Sigma $ when $d \geq 4$) to our present  setting. Denote   $ \{  \pi  _t \}$  the heat semi-group generated by the $L$ in $ \Sigma $ and as above $ \varphi  _0^*$ any positive  eigenfunction of the adjoint elliptic operator $-L^*$ in $ \Sigma $ for the eigenvalue $ \lambda _1$. 
Then we have :

(i) if $d=3$, there is a $t_0>0$ and a function $C:[t_0, \infty ) \to  (1,+ \infty ) $ such that  $\lim _{t \to   \infty }C(t)=1$ and --if $\displaystyle C_0 = (\int _ \Sigma  \varphi  _0(y) \varphi  _0^\ast(y)\, d \sigma  (y))^{-1} $, 
\begin{equation}\label{BanBurEq2} C_0\,C(t)^{-1}\, e^{- \lambda _1t} \, \varphi  _0(x_0)\,  \varphi  _0^{\ast} (y) \leq   \pi  _t(x_0,y) \leq C_0\,C(t)\, e^{- \lambda _1t} \, \varphi  _0(x_0)\,  \varphi  _0^{\ast}(y),\;\; t \geq t_0
\end{equation}
for all   $y \in    \Sigma $. The proof in \cite{bada} Theorem 1 for the Laplacian can be  adapted after one shows that the Cranston-McConnell inequalities \cite{CM} \cite{banCM} (see also  \cite{banProc})  hold for all subdomains $\omega $ of $ \Sigma $: i.e., there is a constant $C=C( \Sigma )$ such that for every $\omega $ and every positive $L$-harmonic function $h$  in $\omega $ one has  $G^\omega (h) \leq C\,  \vert  \omega  \vert   h$.

(ii) for all $d \geq 3$, it is well-known that  (\ref{BanBurEq2}) holds provided $y$ is restricted to a relatively compact subset $A$ of $ \Sigma $ (see  \cite{pincho2} Theorem 1.2 (iii) with  a class of elliptic operators slightly different from ours, see also \cite{pincho}, \cite{pinsky}).

Let us finally also mention that the results in sections \ref{sectjohnreg} and \ref{sectjohncuts} extend  to the present  setting if we restrict to John conditions with $N=1$, where $N$ is the number of poles --recall the needed results in \cite{anc4} require for $N \geq 2$ the symmetry of the underlying operator.  When  $B=0$ and $A$ is symmetric, the restriction $N=1$ can be removed since    $L=L_{ \mathbb R} +\partial ^2_{uu} + b\partial _u$ is then symmetric with respect to the reference measure $\mu  (du,dx)=e^{bu}\, du\, d \sigma _M (x)$ (i.e. $L$ is symmetric in $M\times { \mathbb R}$  equipped with the riemannian metric $g_{(u,x)}(du,dx)=e^{bu}\,( du)^ 2\,   g_M(dx)  $).


\end{document}